\documentclass{amsart}

\usepackage{ricky}

\begin{document}

\title{$p$-adic modular forms of non-integral weight over Shimura curves}

\author{Riccardo Brasca}

\date{October 30, 2011}

\email{\href{mailto:riccardo.brasca@gmail.com}{riccardo.brasca@gmail.com}}
\address{}

\curraddr{Dipartimento di Matematica\\
Universit\`a degli studi di Milano\\
Milan\\
Italy}



\begin{abstract}
In this work, we set up a theory of $p$-adic modular forms over Shimura curves over totally real fields which allows us to consider also non-integral weights. In particular, we define an analogue of the sheaves of $k$-th invariant differentials over the Shimura curves we are interested in, for any $p$-adic character. In this way, we are able to introduce the notion of overconvergent modular form of any $p$-adic weight. Moreover, our sheaves can be put in $p$-adic families over a suitable rigid-analytic space, that parametrizes the weights. Finally, we define Hecke operators, including the $\U$ operator, that acts compactly on the space of overconvergent modular forms. We also construct the eigencurve.
\end{abstract}

\maketitle

\section*{Introduction}
Let $p$ be a prime and let $N>4$ be a fixed positive integer, with $(p,N)=1$. Let $R$ be a separated and complete $\Z_p$-algebra. The first precise definition of the concept of $p$-adic modular form, of level $N$, weight $k \in \Z$, with coefficients in $R$, was given by Serre in \cite{serre_fam}. He identified a modular form with its $q$-expansion, and he defined a $p$-adic modular form as a power series $f \in R[[q]]$ such that there exists a sequence of classical modular forms $\set{f_n}$ that converges $p$-adically to $f$. It turns out that such an $f$ has not always an integral weight, in the sense that $\set{k_n}$, the sequence of weights of $\set{f_n}$, is not eventually constant. To solve this problem, Serre identified an integer $k$ with the map $\Z_p^\ast \to \Z_p^\ast$ that sends $t$ to $t^k$. He then showed that the map $\chi(t) = \varinjlim_n t^{k_n}$ is well defined and moreover $\chi \colon \Z_p^\ast \to \Z_p^\ast$ is a continuous character. This suggests that the weight of a $p$-adic modular form should be a continuous character $\Z_p^\ast \to \Z_p^\ast$. Serre also introduced the notion of an analytic $p$-adic family of modular forms, parametrized by the weight, and showed that the mere existence of the family of the $p$-adic Eisenstein series implies the analyticity of the $p$-adic zeta function. To work with families, it is therefore convenient to introduce the \emph{weight space} $\mc W$. It is a rigid analytic space over $\Q_p$ such that its $K$-points, for any finite extension $K/\Q_p$, are the continuous characters $\Z_p^\ast \to K^\ast$.

Let $Y_1(N)$ be the modular curve of level $N$, over $\Q_p$ ($\Gamma_1(N)$-level structure). Let $X_1(N)$ be the compactification of $Y_1(N)$, we have a universal semi-abelian scheme $\pi \colon A \to X_1(N)$. Let $\underline \omega = \underline \omega_{X_1(N)}$ be the sheaf $e^\ast \Omega^1_{A/X_1(N)}$, where $e \colon X_1(N) \to A$ is the zero section. Classically, a modular form of level $N$ and weight $k$ is defined as a global section of $\underline \omega^{\otimes k}$. In \cite{katz_modu}, Katz gave a geometric interpretation of the notion of a $p$-adic modular form of integral weight (at least in the case $p > 3$, see \cite[Section 2.1]{katz_modu} for what can be done in the case $p=2,3$). For any rational number $0 \leq w < 1$, let $X_1(N)(w)^{\an}$ be the affinoid subdomain of the analytification of $X_1(N)$ defined in \cite{cole_class}, relative to the Eisenstein series $E_{p-1}$. We think about $X_1(N)(w)^{\an}$ as the subset of $X_1(N)^{\an}$ where $E_{p-1}$ has valuation smaller than or equal to 
$w$. The complement of $X_1(N)(0)^{\an}$ is a finite union of discs, called the supersingular discs. Katz introduced the notion of $p$-adic modular form of level $N$, weight $k$, and growth condition $w$: it is a global section of $\underline \omega^{\otimes k}$ on $X_1(N)(w)^{\an}$. A modular form of growth condition $0$ is called a \emph{convergent} modular form, and one of growth condition $w>0$ is called an \emph{overconvergent} modular form. Katz defined also the usual Hecke operators acting on the space of $p$-adic modular forms, including the $\U$ operator, the analogue of the classical $\U_p$ operator of Atkin. Finally, Katz showed that his definition of a $p$-adic modular form generalizes Serre's. In particular, if $f$ is a $p$-adic modular form in the sense of Serre, \emph{of integral weight $k$}, then $f$ can be identified with a convergent $p$-adic modular form, and conversely.

Let $K$ be a finite extension of $\Q_p$. In \cite{serre_ban}, Serre developed Riesz theory for completely continuous endomorphisms of orthonormizable Banach modules over $K$. An example of such a Banach module is provided by the space of $p$-adic modular forms over $K$, of growth condition $w$ and weight any $\chi \colon \Z_p^\ast \to K^\ast$. It is a key fact that, if we consider only modular forms with growth condition $w>0$, the $\U$ operator is completely continuous, so we have a good Riesz theory for it. In \cite{cole_fam}, Coleman developed Riesz theory for a completely continuous operator on a family of orthonormizable Banach modules, generalizing Serre's work. In our work we will need a further generalization of Riesz theory. In \cite{buzz_eigen}, Buzzard showed that the results of Coleman remain true also for Banach modules that are direct summand of an orthonormizable Banach module.

In \cite{cole_fam}, Coleman was able to prove the following theorem, that generalizes \cite{hida} and holds for any $p$.
\begin{teono}
Let $f$ be an overconvergent modular form of weight $k \geq 2$ that is an eigenform for the full Hecke algebra and let $a_p$ be the $\U$-eigenvalue. If $\val(a_p) < k-1$, then $f$ is classical. Furthermore, any such modular form lies in a $p$-adic family of eigenforms over the weight space.
\end{teono}
The first step needed to obtain Coleman's theorem is to define the notion of overconvergent modular form of any weight. A natural approach is to generalize the sheaves $\underline \omega^{\otimes k}$, obtaining the sheaves $\underline \omega^{\otimes \chi}$ on $X_1(w)$, for any $w$ sufficiently small and any $p$-adic weight $\chi$. However, Coleman's approach is completely different. He made a heavy use of the Eisenstein series, in this way he was able to define, and study, the notion of overconvergent modular form of weight $\chi$ through its $q$-expansion. In particular, Coleman did not define the sheaf $\underline \omega^{\otimes \chi}$.

In \cite{over}, Andreatta, Iovita, and Stevens proposed a geometric approach to this problem, as follows. Let $\chi \colon \Z_p^\ast \to K^\ast$ be a continuous character, where $K$ is a finite extension of $\Q_p$ satisfying certain conditions. Then there is a rational number $w>0$ and a locally free sheaf $\Omega_w^\chi$ on $X_1(N)(w)^{\an}$, such that its global sections correspond naturally to $p$-adic modular forms of weight $\chi$ and growth condition $w$, with coefficients in $K$, as defined by Coleman. Furthermore we have Hecke operators and these sheaves can be put in $p$-adic families over the weight space. The same problem is addressed also in \cite{vincent}, where slightly different techniques are used, mainly from Hida theory.

It is natural to try to develop a similar theory for automorphic forms associated to algebraic groups different from $\GL_{2/\Q}$. In this work we study the case of modular forms over certain quaternionic Shimura curves defined starting from a totally real field $F$. These modular forms are related, via the Jacquet-Langlands correspondence, to Hilbert modular forms for $F$, see the introduction of \cite{pay_totally}. The notion of $p$-adic modular form in this context was introduced by Kassaei. The definition is similar to Katz', in particular Kassaei considered only \emph{integral weight}. The goal of this work is to give a geometric definition of quaternionic modular forms of any weight and to prove that these modular forms can be put in families.

Our base ring is $\OP$, the completion of $\mc O_F$ at a prime above $p$. All our objects will be endowed with a natural action of $\OP$. The space of locally analytic character $\OP^\ast \to \C_p^\ast$ has the structure of a rigid analytic variety $\mc W$, called the \emph{weight space}. We fix $K$, a finite extension of $\Fr(\OP)$ satisfying certain technical conditions. Besides the definition of the space of overconvergent modular forms of any weight, our main result is the following (Theorem \ref{teo: eigen}).
\begin{teono}
There is a rigid space $\mc C \subseteq \mc W \times \m A^{1,\rig}_K$, called the eigencurve, such that its $L$-points, where $L$ is a finite extension of $K$, correspond naturally to systems of eigenvalues of overconvergent modular forms defined over $L$. If $x \in \mc C(L)$, let $\mc M_x$ be the set of overconvergent modular forms corresponding to $x$. Then all the elements of $\mc M_x$ have weight $\pi_1(x) \in \mc W(L)$ and the $\U$-operator acts on $\mc M_x$ with eigenvalue $\pi_2(x)^{-1}$.
\end{teono}

In Section~\ref{sec: shim curves}, we define the basic objects of our study, this is due to Carayol in \cite{carayol}. Let $p\neq 2$ denote a fixed rational prime. Let $F$ be a totally real field, with $[F:\Q] > 1$, and let $B$ be a quaternion algebra over $F$ that splits at exactly one infinite place of $F$ and at $\mc P$, a prime of $F$ above $p$. Attached to these data, there is an inverse system of PEL Shimura curves $\set{M_K}$, parametrized by compact open subgroups of $G(\m A^f)$, where $G$ is a reductive algebraic group over $\Q$, defined using $B$. With some assumptions on $K$, we give a precise description of the moduli problem that is solved by $M_K$ over $F_{\mc P}$, the completion of $F$ at $\mc P$, and over $\OP\colonequals  \mc O_{F_{\mc P}}$. The residue field of $\OP$ is denoted with $\kappa$, and we write $q$ for $|\kappa|$.

Section~\ref{sec: hasse} is essentially due to Kassaei. We recall the definition of the analogue of the Hasse invariant in our situation. In this way, we are able to define the space of $p$-adic modular forms of level $\koneH$ (an analogue of level $N$), weight $k \in \Z$, and growth condition $0 \leq w < 1$. In Section~\ref{sec: level Hp}, we recall the theory of the canonical subgroup, as developed in \cite{pay_totally}, and we consider modular forms of higher levels. We can decompose the $p^n$-torsion of the objects of our moduli problems, that are abelian schemes, to define a $p$-divisible group of dimension $1$. In \cite{pay_totally}, this $p$-divisible group is used to define the canonical subgroup. In order to obtain the results we want, we need a generalization of the theory of $p$-divisible group: the theory of $\varpi$-divisible group, where $\varpi$ is a uniformizer of $\OP$. We recall what we need about $\varpi$-divisible groups and we study the $\varpi$-divisible group attached to our abelian 
scheme. Using the canonical subgroup, we are able to show that there is a modular form of level $\koneHpi$ (an analogue of level $Np$), called $E_1$, that is a canonical $q-1$-th root of $E_{q-1}$. This is a new result and it will be essential for our theory. We also obtain some very explicit results about the canonical subgroup, generalizing some results of \cite{coleman_can}.

Section~\ref{sec: dlog} contains the most important technical results of the paper. First of all we show that the trivial analogues of the results of \cite{over} are false in our situation. To be more precise, let $\mc A$ be an object of our moduli problem, of level $\koneHpi$. We have the canonical subgroup $\mc C$ of $\mc A[p]$. In \cite{over}, it is shown that we have a canonical point $\gamma'$ of $\mc C^{\Dual}$ (Cartier dual). One of the most important technical results of \cite{over} is that the image of $\gamma'$ under the map $\dlog$ is congruent, modulo $p^{1-w/(p-1)}$, to $E_1$. Also in our situation we have the canonical point $\gamma'$, but, by Proposition~\ref{prop: gamma' 0}, we have $\dlog(\gamma')=0$ if $\OP$ is sufficiently ramified. This show that we need a different approach. The deep reason for this problem is that all the objects we are interested in are endowed with an action of $\OP$, and we really need to take this action into account. For example, Cartier duality does not work, 
since $\m G_{\operatorname{m}}$ does not have a natural action of $\OP$. What we need is the theory of \emph{group schemes with strict $\OP$-action} as developed by Faltings in \cite{faltings_strict}. Thanks to this theory, we are able to study the $\varpi$-divisible groups attached to our abelian schemes, and we obtain the correct analogue of the conclusions of \cite{over}. In the case $\OP=\Z_p$, we obtain the results of \cite[Sections 2 and 5]{over}.

Having the results of Section~\ref{sec: dlog}, Section~\ref{sec: HT} is similar to \cite{over}. We prove that the homology of the so called Hodge-Tate sequence is killed by a certain power of $\varpi$. This links the Tate module of our abelian schemes to the invariant differentials in a very precise way. Since an elliptic curve admits a canonical principal polarization, all the objects studied in \cite{over} are self-dual. This is not the case in our situation, in particular we need Proposition~\ref{prop: dual}. This lack of self duality makes some of the arguments of Section~\ref{sec: dlog} more delicate than those in \cite{over}.

Section~\ref{sec: sheaf} is the heart of the paper. We prove that we can generalize the definition of the sheaves $\underline \omega^{\otimes k}$ to any locally analytic character $\chi \colon \OP^\ast \to K^\ast$. For any fixed $\chi$, there is a rational $w>0$ and a locally free sheaf $\Omega_w^\chi$ on $\formMoneH(w)^{\rig}$ (this curve is the analogue of $X_1(N)(w)^{\an}$), such that $\Omega_w^\chi = \underline \omega^{\otimes k}$ if $\chi(t)=t^k$. In this way we are able to define the space of quaternionic modular forms of weight $\chi$. In order to define the sheaves $\Omega_w^\chi$, we need to consider the curves $\formMoneHpin{n}(w)^{\rig}$ (analogous of those of level $Np^n$), where $n$ depends on $\chi$. We start by defining a sheaf $\tilde \Omega_w^\chi$ on $\formMoneHpin{n}(w)^{\rig}$. We then show that we have an analogue of the usual diamond operators acting on the push-forward of $\tilde \Omega_w^\chi$ to $\formMoneH(w)^{\rig}$. Taking invariants with respect to these 
operators, we obtain the sheaf $\Omega_w^\chi$. We then consider analytic families. We have an admissible covering $\set{\mc W_r}$ of $\mc W$, made by affinoids. To show that our definition of the sheaves $\Omega_w^\chi$ makes sense, we prove that the $\Omega_w^\chi$'s live in families. More precisely, we prove that there are locally free sheaves $\Omega_{w,r}$ on $\mc W_r \times \formMoneH(w)^{\rig}$, such that $\Omega_w^\chi$ is the pullback of $\Omega_{w,r}$ at the point defined by $\chi$. Furthermore, the $\Omega_{w,r}$'s satisfy various compatibility conditions. This shows that our sheaves really `interpolate' the sheaves $\underline \omega^{\otimes k}$, for various $k$.

In Section~\ref{sec: U} we introduce the Hecke operators $\U$ and $\T_{\mc L}$. These are analogous to the classical $\U$ and $\T_l$ operators. We show that the space of modular forms of weight $\chi$ with coefficients in $K$ is a Banach $K$-module, and that the $\U$ operator is completely continuous when restricted to overconvergent modular forms. We also construct families of the Hecke operators. We prove that a modular form that is an eigenvector for the $\U$ operator (with finite slope) lives in a $p$-adic analytic family of eigenforms. This gives the analogue of the above theorem of Coleman. The main technical results are the following theorems.
\begin{teono}
Let $\chi \colon \OP^\ast \to K^\ast$ be a locally analytic character and assume that $w$ is small enough. Then we have an invertible sheaf $\Omega_w^\chi$ on $\formMoneH(w)^{\rig}_K$ and a completely continuous operator $\U$ on the global sections of $\Omega_w^\chi$. If $\chi(t)=t^k$, then there is a $\U$-equivariant isomorphism between $\Homol^0(\Omega_w^\chi, \formMoneH(w)^{\rig}_K)$ and the space of modular forms of growth condition $w$ as defined in \cite{pay_totally}.
\end{teono}
\begin{teono}
Let $r \geq 0$ be an integer. For any small enough $w$, we have an invertible sheaf $\Omega_{w,r}$ on $\mc W_r \times \formMoneH(w)^{\rig}_K$ such that its pullback to $\formMoneH(w)^{\rig}_K$ at any $\chi \in \mc W_r(K)$ is $\Omega_w^\chi$. We have Hecke operators on $\Omega_{w,r}$. Furthermore, $\U$-eigenforms of finite slope can be deformed.
\end{teono}

As in the classical case, we plan to give a cohomological interpretation of our modular forms. This should allow us to use the powerful language of modular symbols, as for example in \cite{bellaiche}. We finally hope that our approach to define $p$-adic families of overconvergent modular forms can also be used for algebraic groups different from $G$.

\subsection*{Acknowledgements}
I am very grateful to my PhD advisor, Professor Fabrizio Andreatta. I thank him for his patience in the countless conversations we had in his office. I would also like to thank Professor Adrian Iovita that suggested the topic of the paper. I am grateful to Professor Payman Kassaei, who sent me a copy of his PhD thesis when I was not able to understand his papers. I thank Vincent Pilloni for some very useful comments about a preliminary version of this work. His suggestions lead to significant improvements of the paper.

Some of the basic ideas of this work are taken from \cite{over}. I thank Professors Fabrizio Andreatta, Adrian Iovita, and Glenn Stevens for giving me the possibility to read their manuscript. Some of the basic ideas of this work are completely due to them.

\section{Shimura curves and modular forms} \label{sec: shim curves}
In this section we present the basic objects of our work and fix the notations, we will closely follow the presentation given in \cite{pay_totally}, see \cite{carayol} for details.

Let $p\neq 2$ be a prime number, fixed from now on. We fix $F$, a totally real field of degree $N>1$ over $\Q$. We denote with $\tau_1,\ldots,\tau_N$ the various embeddings of $F$ in $\R$ and with $\mc P_1,\ldots,\mc P_m$ the prime ideals of $\mc O_F$ above $p$. The completion of $F$ at $\mc P_i$ will be denoted with $F_{\mc P_i}$. We set $\mc P\colonequals  \mc P_1$, $\tau\colonequals \tau_1$, and $d\colonequals [\FP:\Q_p]$. We write $\OP$ for $\mc O_{\FP}$. Its ramification degree will be denoted with $e$ and its residue degree with $f$. We fix $\varpi$, a uniformizer of $\OP$, and we write $\kappa$ for the residue field $\OP/\varpi\OP$. We write $\val(\cdot)$ for the valuation of $\FP$, normalized by $\val(\varpi)=1$, and we choose $|\cdot|$, an absolute value on $\FP$ compatible with $\val$. We will write $[\cdot] \colon \kappa^\ast \to \OP^\ast$ for the Teichm\"uller character, and we set $[0]=0$. We 
fix $\overline F_{\mc P}$, an algebraic closure of $\FP$, and we 
denote with $\C_p$ its completion. We will use subscripts to denote base change over some base object, that will be clear from the context.

Let $B$ be a quaternion algebra over $F$. We assume that $B$ is split at $\tau$ and at $\mc P$ and that it is ramified at $\tau_2,\ldots,\tau_N$. Let $\lambda < 0$ be a rational number such that $\Q(\sqrt \lambda)$ splits at $p$ and let $E\colonequals F(\sqrt \lambda)$. We embed $E$ in the field of complex numbers via $\tau\colon E \to \C$, where, if $x,y \in F$, $\tau(x+\sqrt\lambda y)=\tau(x)+\sqrt \lambda \tau(y)$. We choose $\mu \in \Q_p$ such that $\mu^2 = \lambda$. We have
\begin{gather} \label{eq: basic deco} 
E \otimes_{\Q}\Q_p \stackrel{\sim}{\longrightarrow} (F_{\mc P_1}\times\cdots\times F_{\mc P_m}) \times (F_{\mc P_1}\times\cdots\times F_{\mc P_m}).
\end{gather}
Composing twice the natural map $E \to E \otimes_{\Q}\Q_p$ with the projection on the first factor, we get a map $E \to \FP$. We use this morphism to define a structure of $E$-algebra on $\FP$.

Let $z\mapsto \bar z$ be the non trivial element of $\Gal(E/F)$. On $D\colonequals B\otimes_F E$, we define an involution $\bar \cdot$, that sends $b \otimes_F z$ to $b' \otimes_F \bar z$, where $\cdot' \colon B \to B$ is the canonical involution of $B$. The underlying $\Q$-vector space of $D$ will be denoted with $V$. We let $D$ act on the left on $V$, by multiplication. We choose $\delta \in D$ such that $\bar \delta = \delta$, and we define an involution $\cdot^\ast \colon D \to D$ by $l^\ast=\delta^{-1}\bar l \delta$. We now take $\alpha \in E$ such that $\bar \alpha = -\alpha$, and we define the symplectic bilinear form
\begin{gather*}
\Theta \colon V \times V \to \Q \\
(v,w) \mapsto \Theta(v,w)= \tr_{E/\Q}(\alpha\tr_{D/E}(v\delta w^\ast))
\end{gather*}

Let $\OB$ be a maximal order of $B$. We fix an isomorphism $\OB \otimes_{\mc O_F} \OP \cong \mat_2(\OP)$. The maximal order of $D$ corresponding to $\OB$ will be denoted with $\OD$, or with $V_{\Z}$ if we want to see it as a lattice in $V$. The isomorphism in~\eqref{eq: basic deco} implies that we have the following isomorphisms
\begin{gather} \label{eq: sec deco}
\begin{array}{cccccccccccc}
\OD \otimes_{\Z} \Z_p & \cong & \mc O_{D_1^1} & \times & \cdots & \times &  \mc O_{D_m^1} & \times & \mc O_{D_1^2} & \times & \cdots & \mc O_{D_m^2}\\
\bigcap & & \bigcap & & & & \bigcap & & \bigcap & & & \bigcap \\
D \otimes_{\Q} \Q_p & \cong & D_1^1 & \times & \cdots & \times &  D_m^1 & \times & D_1^2 & \times & \cdots & D_m^2\\
\end{array}
\end{gather}
We choose $\OD$, $\alpha$, and $\delta$ in such a way that the following conditions are satisfied:
\begin{itemize}
 \item $\OD$ is stable under $l \mapsto l^\ast$;
 \item $\mc O_{D_j^k}$ is a maximal order in $D_j^k$ and $\mc O_{D_1^2}$ is identified with $\mat_2(\OP)$;
 \item $\Theta$ takes integer values on $V_{\Z}$ and it induces a perfect pairing on $V_{\Z_p}\colonequals V_{\Z} \otimes_{\Z} \Z_p$.
\end{itemize}

Let $\mc C$ be a pseudo-abelian category. If $X$ is an object of $\mc C$ with an action of $\OD \otimes_{\Z} \Z_p$, the isomorphism in~\eqref{eq: sec deco} induces a decomposition in $\mc C$
\[
X = X_1^1 \oplus \cdots \oplus X_m^1 \oplus X_1^2 \oplus \cdots X_m^2,
\]
where each $X_j^k$ has an action of $\mc O_{D_j^k}$. Using the idempotents $\left(\begin{array}{cc} 1 & 0 \\ 0 & 0\end{array}\right)$ and $\left(\begin{array}{cc} 0 & 0 \\ 0 & 1\end{array}\right)$ of $\mc O_{D_1^2} \cong \mat_2(\OP)$, we obtain a further decomposition $X_1^2 = X\deco \oplus X_1^{2,2}$. Note that $X\deco$ and $X_1^{2,2}$ are isomorphic. This notation will be used throughout the paper.

Let $G$ be the reductive algebraic group over $\Q$ such that, for any $\Q$-algebra $R$,
\[
G(R)=\set{D\mbox{-linear symplectic similitudes of } (V \otimes_{\Q} R, \Theta\otimes_{\Q} R)}.
\]
We write $\Af$ for the ring of \emph{finite} adele $\Q$. We have
\[
G(\Af) \cong \Q_p^\ast \times \GL_2(\FP) \times (B \otimes_F F_{\mc P_2})^\ast \times \cdots \times (B \otimes_F F_{\mc P_m})^\ast \times G(\m A^{f,p}),
\]
where $\m A^{f,p}$ is the restricted product of the $\Q_l$'s ($l$ prime), with $l \neq p$. We will simply write $\Gamma$ for $(B \otimes_F F_{\mc P_2})^\ast \times \cdots \times (B \otimes_F F_{\mc P_m})^\ast\times G(\m A^{f,p})$. For the rest of the paper, we assume that $K$ is a compact open subgroup of $G(\Af)$ of the form
\[
K=\Z_p^\ast \times K_{\mc P} \times H,
\]
where $K_{\mc P}$ and $H$ are compact open subgroups of $\GL_2(\FP)$ and of $\Gamma$.

We write $\m S$ for the Weil restriction $\Res_{\C/\R}(\m G_{\operatorname{m},\C})$. There is a morphism $h \colon \m S \to G_{\R}$ such that $X$, the $G(\R)$-conjugacy class of $h$, can be identified with $\Hsp$, the Poincar\'e half plane. For $K$ as above, we have the Riemann surface
\[
M_{K}(\C)\colonequals G(\Q)\backslash(G(\Af)\times X)/K.
\]
It can be proved that $M_{K}(\C)$ admits a canonical smooth and proper model over $E$, denoted $M_{K}$. Its base change to $\FP$, denoted again with $M_{K}$, is the Shimura curve we are interested in.

From now on, we assume that $K$ is small enough to keep the lattice $V_{\widehat \Z}\colonequals V_{\Z} \otimes_{\Z} \widehat{\Z} \subseteq V \otimes_{\Q} \Af$ invariant. With this assumption, $M_{K}$ is a fine moduli space, it represents the functor $(\FP\mbox{-algebras})^{\op} \to \mathbf{set}$ that sends $R$, an $\FP$-algebra, to the set of isomorphism classes of quadruples $(A,i,\theta,\bar \alpha)$, where
\begin{enumerate}
 \item $A$ is an abelian scheme over $R$ of relative dimension $4N$, with an action of $\OD$ via $i\colon \OD \to \End_R(A)$ that satisfies:
 \begin{enumerate}
  \item \label{cond: 1a} the projective $R$-module $\lie(A)\deco$ has rank $1$ and $\OP$ acts on it via the natural morphism $\OP \to R$;
  \item for $j\geq 2$, we have $\lie(A)_j^2 =0$;
 \end{enumerate}
 \item $\theta$ is a polarization, of degree prime to $p$, such that the corresponding Rosati involution sends $i(l)$ to $i(l^\ast)$;
 \item $\bar \alpha$ is a $K$-level structure, i.e.\ a class modulo $K$ of a symplectic $\OD$-linear isomorphism $\alpha \colon \widehat{\T}(A) \stackrel{\sim}{\longrightarrow} V_{\widehat{\Z}}$ (locally in the \'etale topology).
\end{enumerate}
Here $\widehat{\T}(A)$ is the product of the Tate modules $\T_l(A)$ ($l$ prime), where each $\T_l(A) = \varprojlim_n A[l^n]$ is considered as an \'etale sheaf. Its symplectic form comes from the Weil pairing composed with $\theta$. Note that $\lie(A)_j^k$ makes sense. We have that $\theta$ induces an isomorphism $A[p^n]_j^1 \stackrel{\sim}{\longrightarrow} (A[p^n]_j^2)^{\Dual}$ (Cartier dual) and that $A[p^n]_j^2$ is \'etale for $j \geq 2$. We write $A[\varpi^n]_1^{2,k}$ for the $\varpi^n$-torsion of $A[p^n]\deco$, and we set $A[\varpi^n]_1^2\colonequals A[\varpi^n]_1^{2,1} \oplus A[\varpi^n]_1^{2,2}$.

In the case $K_{\mc P}$ has some specific form, we can interpret the existence of a $K$-level structure in a more explicit way. We write $\widehat{\rm T}^p(A)$ for $\prod_{l\neq p} \T_l(A)$ and $\widehat{\m Z}^p$ for $\prod_{l \neq p}\Z_l$. We denote $\T_p(A)_2^2 \oplus \cdots \oplus \T_p(A)_m^2$ with ${\rm T}_p^{\mc P}(A)$ and $\left(V_{\Z_p}\right)_2^2 \oplus \cdots \oplus \left(V_{\Z_p}\right)_m^2$ with $W_p^{\mc P}$. Let $\widehat W^p$ be $V_{\Z} \otimes_{\Z} \widehat{\m Z}^p$. We define
\[
\koneH\colonequals \GL_2(\OP),
\]
\[
\kHpin{n}\colonequals \set{\left(\begin{array}{cc} a & b \\ c & d\end{array}\right) \in \GL_2(\OP)\mbox{ s.t. } c\equiv 0 \bmod{\varpi^n}},
\]
and
\[
\koneHpin{n}\colonequals \set{\left(\begin{array}{cc} a & b \\ c & d\end{array}\right) \in \GL_2(\OP)\mbox{ s.t. } a\equiv 1 \bmod{\varpi^n} \mbox{ and } c\equiv 0 \bmod{\varpi^n}}.
\]
In the case $K_{\mc P}=\koneH$, a choice of a level structure is equivalent to a choice of $\bar \alpha^{\mc P}$, where:
\begin{enumerate}
 \item $\bar \alpha^{\mc P}$ is a class, modulo $H$, of $\alpha^{\mc P}=\alpha_p^{\mc P}\oplus \alpha^p$, where $\alpha_p^{\mc P}\colon {\rm T}_p^{\mc P}(A) \stackrel{\sim}{\longrightarrow} W_p^{\mc P}$ is linear and $\alpha^p \colon \widehat{\rm T}^p(A) \stackrel{\sim}{\longrightarrow} \widehat W^p$ is symplectic.
\end{enumerate}
If $K_{\mc P}=\kHpin{n}$, a choice of a level structure is equivalent to a choice of $(C,\bar \alpha^{\mc P})$, where:
\begin{enumerate}
 \item $C$ is a finite and flat subgroup scheme of rank $q^n$ of $A[\varpi^n]\deco$, stable under $\OP$;
 \item $\bar \alpha^{\mc P}$ is as above.
\end{enumerate}
In the case $K_{\mc P}=\koneHpin{n}$, a choice of a level structure is equivalent to a choice of $(Q,\bar \alpha^{\mc P})$, where:
\begin{enumerate}
 \item $Q$ is a point of exact $\OP$-order $\varpi^n$ in $A[\varpi^n]\deco$;
 \item $\bar \alpha^{\mc P}$ is as above.
\end{enumerate}
In these cases, the curves $M_{K}$ will be denoted with $\MoneH$, $\MHpin{n}$, and $\MoneHpin{n}$. They admit canonical proper models over $\OP$, denoted $\intMoneH$, $\intMHpin{n}$, and $\intMoneHpin{n}$. In \cite{carayol} it is proved that $\intMoneH$ is smooth over $\OP$, while the other two curves have semistable reduction.

The curves $\intMoneH$, $\intMHpin{n}$, and $\intMoneHpin{n}$ solve the same moduli problems as the curves $\MoneH$, $\MHpin{n}$, and $\MoneHpin{n}$ do, but now for $\OP$-algebras. The level structure has the same description as above, but now $Q$ is a point of exact $\OP$-order $\varpi^n$ in the sense of Drinfel'd. We have several morphisms between these curves, given by the natural transformations of functors.

The universal objects of the moduli problems of the curves $\MoneH$, $\MHpin{n}$, and $\MoneHpin{n}$ will be denoted with $\AoneH$, $\AHpin{n}$, and $\AoneHpin{n}$. They admit the canonical integral models $\intAoneH$, $\intAHpin{n}$, and $\intAoneHpin{n}$. The morphism $\intAoneH \to \intMoneH$ will be denoted with $\pi$, and its zero section with $e$, we use the same symbols for the other curves.

Let us consider the sheaf $\pi_\ast\Omega^1_{\intAoneH/\intMoneH}$. It has an action of $\OD \otimes_{\Z} \Z_p$, so we can define
\[
\underline \omega \colonequals  \underline \omega_{\koneH} \colonequals  \left(\pi_\ast\Omega^1_{\intAoneH/\intMoneH} \right)\deco.
\]
The definitions of $\underline \omega_{\kHpin{n}}$ and $\underline \omega_{\koneHpin{n}}$ are analogous, we usually drop the subscript. By condition~\eqref{cond: 1a} of the moduli problem, we have that $\underline \omega$ is a locally free sheaf of rank $1$. If $R$ is an $\OP$-algebra, the pullback of $\underline \omega$ to $\Sp(R)$ will be denoted $\underline \omega_R$ or $\underline \omega_{\mc A/R}$, where $\mc A$ is the pullback of the universal object to $\Sp(R)$.
\begin{defi}
Let $R$ be an $\OP$-algebra and $k$ an integer. The space of \emph{modular forms with respect to $D$, level $\koneH$ and weight $k$, with coefficients in $R$}, is defined as
\[
S^D(R,\koneH,k)\colonequals \Homol^0(\intMoneH_R,\underline \omega^{\otimes k}_R).
\]
The definitions of $S^D(R,\kHpin{n},k)$ and $S^D(R,\koneHpin{n},k)$ are similar.
\end{defi}

\section{\texorpdfstring{The Hasse invariant and $\varpi$-adic modular forms}{The Hasse invariant and varpi-adic modular forms}} \label{sec: hasse}
\begin{notation}
We will use the following notation: objects defined over $\OP$ will be denoted using Italics letter, like $\mc A$. The completion along the subscheme defined by $\varpi=0$ will be denoted using the corresponding Gothic letter, like $\mathfrak A$.
\end{notation}

Let $X$ be any $\OP$-scheme (or a formal scheme). Recall that a \emph{$\varpi$-divisible group} $H \to X$ is a Barsotti-Tate group $H$ over $X$, together with an embedding $\OP \hookrightarrow \End(H)$ such that the induced action of $\OP$ on $\lie(H)$ is the one given by $H \to X \to \Sp(\OP)$. If $X$ is connected, there is a unique integer $\height(H)$, called the \emph{height} of $H$, such that $\rk(H[\varpi^n])=q^{n\height(H)}$ for all $n$. Let $\mathfrak X$ be a $\varpi$-adic formal scheme over $\Spf(\OP)$, and let $\mathfrak G \to \mathfrak X$ be a smooth formal group. We say that $\mathfrak G$ is a \emph{formal $\OP$-module} if we have an action of $\OP$ on $\mathfrak G$ such that the action of $\OP$ on $\lie(\mathfrak G)$ is given by $\mathfrak G \to \mathfrak X \to \Spf(\OP)$.

Given $(\mc A,i,\theta,\bar \alpha)$, an object of the moduli problem, with $\mc A$ defined over $R$, we write $\mc A[\varpi^\infty]\deco$ for $\varinjlim_n \mc A[\varpi^n]\deco$, this is the $\varpi$-divisible group associated to $\mc A$. The height of $\mc A[\varpi^\infty]\deco$ is $2$. Let $\mathfrak A$ be the $\varpi$-adic completion of $\mc A$ and let $\widehat{\mc A}$ be the completion of $\mathfrak A$ along its zero section. Then $\widehat{\mc A}\deco$ is a formal $\OP$-module of dimension $1$. If $\widehat{\mc A}\deco$ is a $\varpi$-divisible group, its height $h$ is either $1$ or $2$ and satisfies $q^h = \rk(\widehat{\mc A}\deco [\varpi])$. We say that $(\mc A,i,\theta,\bar \alpha)$, or simply $\mc A$, is \emph{ordinary} if $\widehat{\mc A} \deco$ has height $1$. If $\widehat{\mc A} \deco$ has height $2$, we say that $(\mc A,i,\theta,\bar \alpha)$ is \emph{supersingular}.

With the above notations, suppose that $\widehat{\mc A}\deco$ is coordinatizable. It is proved in \cite[Proposition~4.3]{pay_totally}, that we can find a coordinate $x_R$ on $\widehat{\mc A}\deco$ such that the action of $\varpi$ has the form
\[
[\varpi](x_R)=\varpi x_R + a_R x_R^q + \sum_{j=2}^{\infty}c_jx_R^{j(q-1)+1},
\]
where $a$, $c_j$ are in $R$ and $c_j \in \varpi R$ unless $j \equiv 1 \bmod{q}$. If we assume that $\varpi=0$ in $R$, the various $a_R$ glue together to define $\mathbf{H}$, a modular form of level $\koneH$ and weight $q-1$, defined over $\kappa$, that is called the \emph{Hasse invariant}. If $W=\Sp(R)$ is an open affine of $\intMoneH_{\kappa}$ and we denote with $\omega$ the differential dual to the coordinate $x_R$ defined above, we have $\mathbf{H}_{|W}=a_R \omega^{\otimes q-1}$.

We assume $q>3$ (but see Remark \ref{rmk: q}), so by \cite[Proposition 7.2]{pay_totally}, the Hasse invariant can be lifted to a modular form of level $\koneH$ and weight $q-1$, defined over $\OP$. We choose such a lifting, called $E_{q-1}$: in \cite[Corollary 13.2]{pay_totally}, it is shown that all the theory does not depend on this choice. Over $\Sp(R)$, we can write ${E_{q-1}}_{|\Sp(R)}=E\omega^{\otimes q-1}$, with $E \in R$. By \cite[Proposition~6.2]{pay_totally}, we have $a_R \equiv E \bmod{\varpi}$.
\begin{prop}[{(\cite[Proposition~6.1]{pay_totally})}] \label{prop: H vanish ok}
Let $R$ be a $\kappa$-algebra and let $(\mc A,i,\theta,\bar \alpha)$ be an object of the moduli problem, with $\mc A$ defined over $R$ and let $z$ be a geometric point of $\Sp(R)$. Then the pullback of $\mathbf{H}$ to $\intMoneH_R$ vanishes at $z$ if and only if the pullback of $\mc A$ to $z$ is supersingular.
\end{prop}
We now move on to $\varpi$-adic modular form. Let $V$ be a finite extension of $\OP$ and let $w$ a rational number such that there is an element of $V$, denoted $\varpi^w$, of valuation $w$. We define
\[
\intMoneH(w)_V\colonequals \Sp_{\intMoneH_V}(\Sym(\underline \omega^{\otimes q-1}_V)/(E_{q-1}-\varpi^w)).
\]
\begin{rmk} \label{rmk: q}
If $q=2,3$, we can still define $\formMoneH(w)$, that is the $\varpi$-adic completion of $\intMoneH(w)$, using the notion of `measure of singularity' introduced in \cite[Section~2]{pay_curve}. Over $\Spf(R)$, we can take for $E$ any element of $R$ such that $E \omega^{\otimes q-1}$ is a generator of $\underline \omega ^{\otimes q-1}$. It follows that all `local' results of the paper, i.e.\ those that do not involve the existence of $E_{q-1}$, remain true also if $q=2,3$. We leave the details to the interested reader.
\end{rmk}
\begin{defi}
Let $V$ and $w$ as above. The space of $\varpi$-adic modular forms with respect to $D$, level $\koneH$, weight $k$ and growth condition $w$, with coefficients in $V$, is defined as
\[
S^D(V,w,\koneH,k)\colonequals \Homol^0(\formMoneH(w)_V,\underline \omega^{\otimes k}).
\]
\end{defi}
Let $K$ be the fraction field of $V$. The rigidification of the map $\formMoneH(w)_V\to \formMoneH_V$ is the immersion $\formMoneH_V^{\rig}(w) \hookrightarrow \formMoneH_V^{\rig}$, where $\formMoneH_V^{\rig}(w)$ is the affinoid subdomain of $\formMoneH_V^{\rig}$ defined by Coleman in \cite{cole_class}, relative to $E_{q-1}$. We call $\formMoneH_V(0)^{\rig}$ the \emph{ordinary locus}, it is an affinoid subdomain of $\formMoneH_V^{\rig}$: its complement is a finite union of the \emph{supersingular discs}.

By rigid GAGA, elements of $S^D(V, \koneH, k)_K$ correspond to sections of $\underline \omega^{\otimes k}$ over $\formMoneH_V^{\rig}$, while elements of $S^D(V,w,\koneH,k)_K$ correspond to sections over $\formMoneH_V(w)^{\rig}$. Elements of $S^D(V,0,\koneH,k)_K$ are called \emph{convergent} modular with coefficients in $K$, while the elements of $S^D(V,w,\koneH,k)_K$, for $w>0$, are called \emph{overconvergent} modular forms.

\section{\texorpdfstring{The canonical subgroup and modular forms of level $\koneHpi$}{The canonical subgroup and modular forms of level K(Hvarpi)}} \label{sec: level Hp}
\begin{notation}
Let $V$, $K$, and $w$ be as above. From now on, we will work over $V$, so we will consider the base change to $V$, or to $K$, of the various objects defined so far. For simplicity we will omit the subscripts $_V$ and $_K$.
\end{notation}
From now on we assume that $0 \leq w < \frac{q}{q+1}$. By \cite[Theorem~10.1]{pay_totally}, the $q$-torsion of any $\mc A$ as above admits a canonical subgroup, that we call $\mc C$. We have that $\mc C\deco$ is killed by $\varpi$. Since $\mc C\deco$ has order $q$, we can use it to define a morphism $\formMoneH(w) \to \formMHpi$. Its rigidification is a section, that is defined over $\formMoneH(w)^{\rig}$, of the morphism $\formMHpi^{\rig} \to \formMoneH^{\rig}$. We define $\formMoneHpi(w)^{\rig}$ as the inverse image of $\formMoneH(w)^{\rig}$ with respect to the map $\formMoneHpi^{\rig} \to \formMHpi^{\rig}$. It is an affinoid subdomain of $\formMoneHpi^{\rig}$ with a finite and \'etale map to $\formMoneH(w)^{\rig}$.

We assume that $V$ contains an element, denoted $(-\varpi)^{1/(q-1)}$, whose $q-1$-th power is $-\varpi$. Let $\mathfrak U=\Spf(R)$ be an open affine of $\formMoneH(w)$. We write $\mathfrak U^{\rig}=\Spm(R_K)$ for its rigid analytic fiber. Since the morphism $\formMoneHpi(w)^{\rig} \to \formMoneH(w)^{\rig}$ is finite and \'etale, the inverse image of $\mathfrak U^{\rig}$ is an affinoid, $\mathfrak V^{\rig}=\Spm(S_K)$, with $R_K \to S_K$ finite and \'etale. Let $S$ be the normalization of $R$ in $S_K$, and let $\mathfrak V$ be $\Spf(S)$. Note that $S$ is $\varpi$-adically complete. The various $\mathfrak V$'s glue to define a formal scheme $\formMoneHpi(w)$, with a morphism to $\formMoneH(w)$. We have the following
\begin{lemma}
The rigid analytic fiber of $\formMoneHpi(w)$ is $\formMoneHpi(w)^{\rig}$. Furthermore, the rigidification of $\formMoneHpi(w) \to \formMoneH(w)$ is the map $\formMoneHpi(w)^{\rig} \to \formMoneH(w)^{\rig}$ defined above.
\end{lemma}
By definition, $\formMoneHpi(w)$ is the normalization of $\formMoneH(w)$ in $\formMoneHpi(w)^{\rig}$, that is a finite extension of its generic fiber.
\begin{prop} \label{prop: mod Hpi}
Let $S$ be a normal and $\varpi$-adically complete $V$-algebra. There is a natural bijection between $\formMoneHpi(w)(S)$ and the set of isomorphism classes of quintuples $(\mc A, i, \theta, \bar\alpha, Y)$, where:
\begin{itemize}
 \item $(\mc A, i, \theta, \bar\alpha)$ is an object of the moduli problem, with $\mc A$ defined over $S$, of $\intMoneHpi$ and the canonical $S$-point of $\mc A[\varpi]\deco$ generate, as $\OP$-module, the canonical subgroup of $\mc A[\varpi]$;
 \item $Y$ is a section of $\underline \omega_{\mc A/S}^{\otimes 1-q}$ that satisfies $YE_{q-1}=\varpi^w$.
\end{itemize}
\end{prop}
\begin{proof}
This is proved in exactly the same way as \cite[Lemma~3.1]{over}.
\end{proof}
\begin{defi}
We define the space of $\varpi$-adic modular forms with respect to $D$, level $\koneHpi$, weight $k$ and growth condition $w$, with coefficients in $V$, as
\[
S^D(V,w,\koneHpi,k)\colonequals \Homol^0(\formMoneHpi(w),\underline \omega^{\otimes k}).
\]
\end{defi}
Note that we have $S^D(V,w,\koneHpi,k)_K = \Homol^0(\formMoneHpi(w)^{\rig},\underline \omega^{\otimes k})$. We have a natural map $S^D(V,w,\koneH,k) \to S^D(V,w,\koneHpi,k)$. The image of $E_{q-1}$ will be denoted with the same symbol. We fix an open affine $\Spf(R)$ of $\formMoneH(w)$, with associated abelian scheme $\mc A$. Let $x$ be a coordinate on $\widehat{\mc A}\deco$ as in Section~\ref{sec: hasse}. We write ${E_{q-1}}_{|\Spf(R)}=E \omega^{\otimes q-1}$. Let $\Spf(S)$ be the base change of $\Spf(R)$ to $\formMoneHpi(w)$. We need to briefly recall how $\mc C\deco$ is constructed, see \cite[Lemma~10.2]{pay_totally} for details. There is $b \in R$ such that $E = a + b\varpi$ and we can write $(a+\varpi b)y=\varpi^w$, with $y \in R$. We set $r_1\colonequals -\varpi/\varpi^w \in V$ and $t_0\colonequals  r_1y/(1+r_1by) \in R$. We have that $\mc C\deco$, as a scheme, is $\Sp(R[[x]]/(x^q-t_{\operatorname{can}}x))$, where $t_{\operatorname{can}}=t_0(1-t_\infty)$. Here $t_\infty$ is an element of $r_2R$, where 
$r_2 \in V$ has positive valuation. Since $t_{\operatorname{can}}$ is topologically nilpotent, we have an isomorphism $\mc C \deco \cong \Sp(R[x]/(x^q-t_{\operatorname{can}}x))$. It follows that there is $r \in R$, such that $x \mapsto r x$ gives
\[
\mc C\deco \cong \Sp(R[x]/(x^q+\frac{\varpi}{E}x)).
\]
\begin{prop} \label{prop: E exist}
There is $E_1 \in S^D(V,w,\koneHpi,1)$ such that $E_1^{q-1}=E_{q-1}$.
\end{prop}
\begin{proof}
By Proposition~\ref{prop: mod Hpi}, the equation $x^{q-1}+\frac{\varpi}{E}=0$ has a canonical solution $\alpha \in S$. Consider the element $E^{1/(q-1)}\colonequals \frac{(-\varpi)^{1/(q-1)}}{\alpha} \in S_K$: it is a canonical $q-1$-th root of $E$ in $S_K$, that lies in $S$ by normality. For the various $R$'s, these roots glue to define the required modular form.
\end{proof}
\subsection{Raynaud theory} \label{subsec: raynaud}
We will continue to work locally for all this section, using the notations introduced above. In this section we find it convenient to denote the Teichm\"uller character $[\cdot]$ with $\chi_1(\cdot)$ (see below).

Let $M$ be the set of multiplicative characters $\chi \colon \kappa^\ast \to \OP^\ast$, extended to the whole $\kappa$ by $\chi(0)=0$. Following Raynaud in \cite{raynaud_type}, we say that $\chi \in M$ is a \emph{fundamental character} if the map $\kappa \to \OP \to \OP/\varpi \OP=\kappa$ is a field homomorphism. If $\chi$ satisfies this condition, all fundamental characters are of the form $z \mapsto \chi(z)^{p^i}$. We denote this latter character with $\chi_{p^i}$, where $i \in \Z/f\Z$ (in \cite{raynaud_type}, $\chi_{p^i}$ is denoted with $\chi_i$). Furthermore we can assume that $\chi_{p^{i+1}}=\chi_{p^i}^p$ and that $\chi_1$ is the Teichm\"uller character. Any $\chi \in M$ can be decomposed as $\chi = \prod_{i \in \Z/f\Z} \chi_{p^i}^{n_i}$, with $0 \leq n_i \leq p-1$. If $\chi \neq 1$, this decomposition is unique. In this way we get a bijection between $M\setminus \set{1}$ and $\set{1, \ldots, q-1}$ given by $\chi = \prod_i \chi_{p^i}^{n_i} \mapsto \sum_i n_i p^i$. We write $\chi_i$ for the inverse 
image of $i$ (this notation is different from the one in \cite{raynaud_type}).

We have shown above that $\mc C \deco$ is, as a scheme, $\Sp(R[x]/(x^q+\frac{\varpi}{E}x))$. Let $w$ and $w_i \colonequals w_{\chi_i}$ be the universal constants introduced by Raynaud. By \cite[Corollary 1.5.1]{raynaud_type}, there is $(\gamma_i,\delta_i)_{i \in \Z/f\Z} \in R^{\Z/fZ}$ such that $\gamma_i\delta_i = w$ and
\[
\mc C \deco \cong \Sp(R[x_i]_{i \in \Z/f\Z}/(x_i^p - \delta_i x_{i+1})).
\]
Since the action of $\OP$ on $\mc C \deco$ is strict, we have $\delta_i \in R^\ast$ for all $i \not \equiv f-1$. We can thus write $x_i = \delta_{i-1}^{-1}\cdots \delta_0^{-p^{i-1}} x_0^{p^i}$ for all $i \not \equiv 0$. It follows that we can assume $x=x_0$ and $-\frac{\varpi}{E}=\delta_0^{p^{f-1}}\cdots \delta_{f-1}$. Furthermore, the action of $z \in \kappa$ on $R[x]/(x^q+\frac{\varpi}{E}x)$ is given by $x \mapsto \chi_1(z)x$ and the module of invariant differentials of $\mc C \deco$ is isomorphic to $R/\frac{\varpi}{E}R \di(x)$ (see \cite[Section 1.1.2]{fargues_hodge}).

Let $h_i$ be the smallest integer, with $0 < h_i \leq f$, such that $p^{f-h_i}$ divides $i$. By \cite[Proposition~1.3.1]{raynaud_type}, there is a unit $u \in \OP^\ast$ such that $w=pu=\varpi^eu$. A calculation with \cite[Corollary 1.5.1]{raynaud_type}, gives the following
\begin{prop} \label{prop: comultiplication}
The comultiplication in $\mc C\deco$ is given by
\begin{gather*}
c(x)=x \otimes 1 + 1 \otimes x - \varpi^{e-1}uE\sum_{i=1}^{q-1}\frac{w^{h_i-1}}{w_i w_{q-i}} x^i\otimes x^{q-i}.
\end{gather*}
\end{prop}
\begin{rmk} \label{rmk: relations OT}
Do not confuse our $w_i$'s with Raynaud's ones, that are all equal and are denoted with $w$ here. Since $\chi_{p^i}$ is a fundamental character for each $i$, we have $w_{p^i}=1$ for $i=0,\ldots,f-1$. Suppose that $\varpi = p$ and $f = 1$. Let us write $w'_i$ for the universal constants used in \cite{coleman_can}. We have $w_i = (-1)^{i+1}\frac{w'_i}{(p-1)^{i-1}}$ and $w=\frac{w'_p}{(p-1)^{p-1}}$, so our description of the comultiplication is exactly the same as the one given in \cite{coleman_can} for the canonical subgroup of an elliptic curve.
\end{rmk}
We assume for the rest of the section that $V$ contains $\zeta_p$, a fixed primitive $p$-th root of unity. Since the base change of $\mc C \deco$ to $S_K$ is a constant group scheme, with associated abstract group $\kappa$, the Hopf algebra of $(\mc C\deco)_{S_K}$ is isomorphic to the algebra of $S_K$-valued functions on $\kappa$. Let $\varepsilon_z$, for $z \in \kappa$, denote the characteristic function of $\set{z}
$. We have a natural isomorphism $S_K[x]/(x^q+\frac{\varpi}{E}x) \to \bigoplus_{z \in \kappa} S_K \varepsilon_z$ sending $x$ to $\sum_{z \in \kappa}\chi_1(z)\alpha\varepsilon_z$. Here $\alpha$ is the root of $-\frac{\varpi}{E}$ given in the proof of Proposition~\ref{prop: E exist}. The Hopf algebra of $(\mc C\deco)_{S_K}^{\Dual}$ is isomorphic to $S_K[\kappa]$, the group algebra of $\kappa$ with coefficients in $S_K$. The canonical base of $S_K[\kappa]$ will be denoted $\set{\mathbf z}_{z \in \kappa}$. Using $\zeta_p$, we can identify $\mathds F_p$ with $\mu_p(V)$. In particular, the trace map $\tr_{\kappa/\mathds F_p}$ can be seen as a morphism $\Psi \colon \kappa \to \mu_p(V)$. We obtain a morphism of group schemes $(\mc C\deco)_{S_K} \to \mu_p$. This morphism corresponds to the $S_K$-point of $(\mc C\deco)_{S_K}^{\Dual}$ given by $\mathbf z \mapsto \Psi(z)$, and comes from
\begin{align*}
\eta \colon S_K[y]/(y^p-1) &\to \bigoplus_{z \in \kappa} S_K \varepsilon_z \cong S_K[x]/(x^q+\frac{\varpi}{E}x) \\
y &\mapsto \sum_{z \in \kappa}\Psi(z)\varepsilon_z
\end{align*}
Let $e_{\chi_i}$, with $0 < i < q-1$, be $\sum_{z \in \kappa}\chi^{-1}_i(z)\mathbf z$ and let $e_{\chi_{q-1}}$ be $\sum_{z \in \kappa}\mathbf z - q \mathbf 0$. We have
\[
\sum_{z \in \kappa}\Psi(z)\varepsilon_z=\varepsilon_0+\frac{1}{q-1}\sum_{i=1}^{q-1}e_{\chi_i}\left(\sum_{z \in \kappa^\ast}\Psi(z)\varepsilon_z \right)\varepsilon_{\chi_i}= 1+\frac{1}{q-1}\sum_{i=1}^{q-1}g(\chi_i)\alpha^{-i}x^i,
\]
where $g(\chi_i)$ is the Gauss sum associated to $\chi_i^{-1}$ and $\Psi$.

If $0 \leq i \leq q-1$ is an integer, written in base $p$ as $i=\sum_{k=0}^{f-1}i_kp^k$, we define $s(i)$ to be $i_0+\ldots + i_{f-1}$. If $i \neq 0$, by \cite[page~251]{raynaud_type}, we have
\[
w_i=w_{\chi_i}=w_{\underbrace{\chi_1,\ldots,\chi_1}_{i_0 \text{ times}},\ldots,\underbrace{\chi_{p^{f-1}},\ldots,\chi_{p^{f-1}}}_{i_{f-1} \text{ times}}}=\frac{g(\chi_1)^{i_0}\cdots g(\chi_{p^{f-1}})^{i_{f-1}}}{(q-1)^{s(i)-1}g(\chi_i)}.
\]
Since $g(\chi_{p^k})=g(\chi_1)$ for every $k$, we obtain
\begin{gather} \label{eq: g}
g(\chi_i)=\frac{1}{(q-1)^{s(i)-1}}\frac{g(\chi_1)^{s(i)}}{w_i}.
\end{gather}
For $k=0,\ldots,f-1$, let $\beta_k\colonequals g(\chi_1)\alpha^{-p^k}$. Writing $i=\sum_{k=0}^{f-1}i_kp^k$, we have
\[
\eta(y)= 1 + \sum_{i=1}^{q-1}\frac{1}{(q-1)^{s(i)}}\frac{x^i}{w_i}\prod_{k=0}^{f-1}\beta_k^{i_k}.
\]
\begin{prop} \label{prop: gamma exists}
The morphism $\eta \colon R[y]/(y^p-1) \to R[x]/(x^q+\frac{\varpi}{E}x)$ induces a canonical $S$-point of $(\mc C \deco)^{\Dual}$. Its base change to $S_K$, denoted with $\gamma'$, is a $\kappa$-generator of $(\mc C \deco)^{\Dual}(S_K)$.
\end{prop}
\begin{proof}
It is enough to show that each $\beta_k$ is in $S$. By equation~\eqref{eq: g}, the valuation of $g(\chi_1)$ is $\frac{e}{p-1}$, so in $S_K$ we can write $\beta_k = v\varpi^{\frac{e}{p-1}-\frac{p^k}{q-1}}(\varpi^{\frac{1}{q-1}}\alpha^{-1})^{p^k}$, where $v$ is a unit of $V$. The claim follows.
\end{proof}
\begin{rmk}
Using the relations between our $w_i$'s and the universal constants used by Coleman given in Remark~\ref{rmk: relations OT}, we see that, in the case $f=1$ and $\varpi=p$, our morphism is exactly the one defined in \cite[Proposition~5.2]{over}.
\end{rmk}
\begin{prop} \label{prop: descr diff}
Let $h \colon \mc C \deco \to \widehat{\mc A}[\varpi]\deco$ be the natural map. In $\Omega^1_{\mc C \deco /R}$, we have the equality
\[
h^\ast(\omega) = \frac{\di(x)}{1-\varpi^{e-1}uE\frac{w^{f-1}}{w_{q-1}}x^{q-1}}.
\]
Furthermore, if we write $\underline \omega_{\mc C \deco/R} \cong R/\frac{\varpi}{E}R \di(x)$, we have $h^\ast(\omega)=\di(x)$.
\end{prop}
\begin{proof}
It is convenient to write the comultiplication $c(x)$ of $\mc C \deco$ as $F(X,Y)$, where $X= x \otimes 1$ and $Y=1\otimes x$. Let $f(X) \di(X)$ be an invariant differential, we have
\[
f(X) \di(X) + f(Y) \di (Y) = f(F(X,Y))(\frac{\partial}{\partial X}F(X,Y)+\frac{\partial}{\partial Y}F(X,Y)),
\]
so, comparing the coefficients of $\di(Y)$ in the two sides of the equation and setting $Y=0$, we find that $f(0) \equiv f(X)(1-QX^{q-1}) \bmod{\frac{\varpi}{E}}$, where $Q\colonequals \varpi^{e-1}uE\frac{w^{f-1}}{w_{q-1}}$. By \cite[page~251]{raynaud_type}, $w=\frac{g(\chi_1)^{p-1}}{(q-1)^{p-1}}$, so we obtain that $(1-(q-1)QX^{q-1})(1-QX^{q-1})=1$, so any invariant differential on $\mc C \deco$ has the form
\[
\frac{r\di(x)}{1-Qx^{q-1}},
\]
for some $r \in R/\frac{\varpi}{E}R$. Since $\omega$ is a differential dual to $x$, we have $\omega= f(x) \di(x)$, with $f \equiv 1 \bmod{x}$ and the first part of the proposition follows. The last statement is a consequence of the fact that the counit of $\mc C \deco \cong \Sp(R[x]/(x^q+\frac{\varpi}{E}x))$ is the map $x \mapsto 0$.
\end{proof}
\begin{rmk}
If $\varpi=p$ and $f=1$, we have $\omega = \frac{\di(x)}{1+\frac{E}{p-1}}$, but $p\di(x)=0$, so $h^\ast(\omega) = \frac{\di (x)}{1-E}$. 
\end{rmk}

\section{\texorpdfstring{The map $\dlog$}{The map dlog}} \label{sec: dlog}
Let $\set{\Spf(R_i)_{i \in I}}$ be a covering of $\formMoneH(w)$ by small affine formal schemes (in the sense of \cite{brinon}). Our local situation will be the following: we choose one of the $R_i$'s, called simply $R$. Its pullback to $\formMoneHpi(w)$ will be denoted with $\Spf(S)$. We assume that $\underline \omega_{\mc A / R}= \left(\pi_\ast \Omega^1_{\mc A / R}\right)\deco$ is a free $R$-module, generated by $\omega$, and we write ${E_{q-1}}_{|\Spf(R)}=E \omega^{\otimes q-1}$. Let $\eta = \Sp(\m K)$ be a generic geometric point of $\Sp(R)$, we write $\mc G$ for $\pi_1(\Sp(R_K), \eta_K)$. We denote with $\overline R$ the direct limit of all $R$-algebras $T \subseteq \m K$ which are normal and such that $T_K$ is finite and \'etale over $R_K$. Let $\widehat{\overline R}$ be the $\varpi$-adic completion of $\overline R$. We are going to use the map $\dlog$, for the definition and basic properties look at \cite[Section~2]{over}, or see below.
\begin{prop} \label{prop: gamma' 0}
If $e$ is big enough, then $\gamma'$ is in the kernel of the map
\[
\dlog \colon (\mc C \deco)^{\Dual}(S_K) \to \underline \omega_{(\mc C \deco)_S / S} \otimes_S S/pS.
\]
\end{prop}
\begin{proof}
We have $\dlog(\gamma')=d\eta(y)/\eta(y)$. By the explicit definition of $\gamma'$, it is enough to prove that $\varpi$ divides $\beta_k$ for each $k=0,\ldots, f-1$ (see Section~\ref{subsec: raynaud} for the definition of $\beta_k$ and $\eta$). In the proof of Proposition~\ref{prop: gamma exists} we have shown that $\varpi^{\frac{e}{p-1}-\frac{p^k}{q-1}}$ divides $\beta_k$. If $e$ is big enough, this implies that $\varpi$ divides $\beta_k$ as required.
\end{proof}
\begin{rmk} \label{rmk: over false}
The above proposition shows that, in general, the analogue of \cite[Proposition~5.2]{over} is not true in our situation. Since we work with $\varpi$-divisible groups rather than with $p$-divisible groups, to obtain results similar to those of \cite{over}, we need a theory that takes into account the action of $\OP$. For example, in Proposition~\ref{prop: gamma' 0}, we use that $\mc C \deco$ is killed by $p$, but we even know that $\varpi \mc C \deco = 0$. Note that there is no action of $\OP$ on $\m G_{\operatorname{m},S}$, so Cartier duality does not suffice. We need the theory of \emph{group schemes with strict $\OP$-action} as developed in \cite{faltings_strict}.
\end{rmk}
On the power series ring $R[[x]]$ there is a unique action of $\OP$ such that the multiplication by $\varpi$ has the form $[\varpi](x)=x^q+\varpi x$ and the action on the Lie algebra is the one induced by $\OP \to R$. This is the so called \emph{Lubin-Tate $\varpi$-divisible group}, denoted $\mc{LT}$. The action of $\OP$ on $R[[x]]/(x^q+\varpi x)$ factors through $\kappa$, and $z \in \kappa$ sends $x$ to $[z]x$. Let $G$ be a finite and flat group scheme over $R$ with an action of $\OP$ (we will always assume the condition on the action on the Lie algebra), and let us suppose that $G$ is killed by $\varpi^n$ for some integer $n$. The functor $(R \mbox{-algebras})^{\op} \to \mathbf{grp}$, that sends $T$ to $\hom_{\OP}(G_T,\mc{LT}_T)$, is representable by a finite and flat group scheme over $R$, with an action of $\OP$, denoted $G^\vee$, see \cite[Sections 3 and 5]{faltings_strict}.

Let $G$ be a $\varpi$-divisible group and let $H$ be a sub $\OP$-module of $\T_\varpi(G^\vee)\colonequals \varprojlim_n G^\vee[\varpi^n](\overline R_K)$ (this is the \emph{Tate module} of $G^\vee$). By duality between $G$ and $G^\vee$, we obtain $H^\perp$, the orthogonal of $H$, that is a sub $\OP$-module of $\T_\varpi(G)$.

If $D \subseteq G[\varpi^n](\overline R_K)$ is a sub $\OP$-module, we write $D^{\cl}$ for the schematic closure of $D$ in $G[\varpi^n]$. Let $R$ be a discrete valuation ring, whose valuation extends the one of $\OP$, so $D^{\cl}$ and $(D^\perp)^{\cl}$ are group schemes. By \cite[Proposition~1]{fargues_hodge}, we have $(D^{\cl})^\vee \cong G[\varpi^n]^\vee/(D^\perp)^{\cl}$.

Let $W$ be a normal Noetherian $R$-algebra, without $\varpi$-torsion. Let $G$ be a group scheme with an action of $\OP$, and let $\underline \omega_{G/R}$ be the module of invariant differential of $G$. If $G$ is killed by $\varpi^n$, we define a map
\[
\dlog_G \colonequals  \dlog_{G,W} \colon G^\vee(W_K) \to \underline \omega_{G/R} \otimes_R W/\varpi^n W
\]
in the following way: given $x$, a $W_K$-valued point of $G^\vee$, it extends to a $W$-valued point of $G^\vee$, called again $x$. It gives a group scheme homomorphism (that respects the action of $\OP$) $f_x \colon G \to \mc{LT}$. We define
\[
\dlog_{G,W}(x) \colonequals  f_x^\ast \di(T).
\]
The map $\dlog$ satisfies various functoriality properties, see \cite[Lemma 2.1]{over}.

We can take $G=\mc A[\varpi^n] \deco$, and we obtain the map
\[
\dlog_{n,W} \colon (\mc A[\varpi^n]\deco)^\vee(W_K) \to \underline \omega_{\mc A[\varpi^n]\deco} \otimes_R W/\varpi^n W.
\]
Taking the direct limit over all $W$ as above, we get the map
\[
\dlog_{n, \mc A} \colon (\mc A[\varpi^n]\deco)^\vee (\overline R_K) \to \underline \omega_{\mc A/R} \otimes_R \overline R/\varpi^n \overline R.
\]
Finally, taking the projective limit, we obtain the map
\[
\dlog_{\mc A} \colon \T_\varpi((\mc A[\varpi^\infty]\deco)^\vee) \to \underline \omega_{\mc A/R} \otimes_R \widehat{\overline R}.
\]
Suppose that $R$ is a discrete valuation ring, whose valuation extends the one of $\OP$. From $\dlog_{\mc A}$, we obtain the maps $\dlog_{n, \widehat{\mc A}}$ and the map
\[
\dlog_{\widehat{\mc A}} \colon \T_{\varpi}((\widehat{\mc A}\deco)^\vee) \to \underline \omega_{\mc A/R} \otimes_R \widehat{\overline R},
\]
\begin{rmk} \label{rmk: newton}
Let us assume that $\mc A$ is supersingular. We have that the Newton polygon of $[\varpi](x)$ is the convex hull of the points
\[
(0, +\infty), \; (1,1), \; (q,\val(a)), \mbox{ and } (q^2,0).
\]
In particular, we see that the roots of $[\varpi](x)$ corresponding to points of $\mc C \deco$ are those with biggest valuation, as in \cite{fargues_hodge}.
\end{rmk}
\begin{prop} \label{prop: ker dlog}
Let us suppose that $w \leq \frac{1}{q}$. Taking the quotient over the kernel, the map $\dlog_{1,\mc A}$ factors through a map, denoted again
\[
\dlog_{1,\mc A} \colon (\mc C \deco)^\vee(\overline R_K) \to \underline \omega_{\mc A/R}\otimes_R \overline R/\varpi \overline R.
\]
Furthermore, the cokernel of the base change to $\overline R/\varpi \overline R$, over $\kappa$, of this map is killed by $\varpi^v$, where $v\colonequals \frac{w}{q-1}$. In particular we have $\ker(\dlog_{1,\mc A})= (\mc C \deco(\overline R_K))^\perp$.
\end{prop}
\begin{proof}
As in the proof of \cite[Proposition~5.1]{over}, we can assume that $R$ is a discrete valuation ring, whose valuation extends the one of $V$. We can prove the proposition with $\mc A\deco$ replaced by $\widehat{\mc A}\deco$. Indeed, if $\mc A$ is supersingular we have $\widehat{\mc A}[\varpi]\deco = \mc A[\varpi]\deco$, while in the ordinary case we can use \cite[Lemma~1]{fargues_hodge}.

First of all we prove the proposition in the case $\widehat{\mc A}\deco$ has height 2. Let $y \in \widehat{\mc A}\deco[\varpi]^\vee(\overline R_K)$ and let $\mc D \subseteq \widehat{\mc A}\deco[\varpi]^\vee(\overline R_K)$ be the $\OP$-module generated by $y$. We have $\underline \omega_{(\mc D^{\cl})^\vee} \cong R/\gamma R$, with $\val(\gamma)=1-\sum \val(z)$, where the sum is over $\mc D^\perp \setminus \set{0}$ (see \cite{fargues_hodge}). Since the map
\[
\overline R/\gamma \overline R \cong\underline \omega_{(\mc D^{\cl})^\vee /R} \otimes_R \overline R / \varpi \overline R \hookrightarrow \underline \omega_{\mc A/R} \otimes_R \overline R/\varpi \overline R \cong \overline R /\varpi \overline R
\]
is the multiplication by $\frac{\varpi}{\gamma}$, it is injective. In particular, we have a commutative diagram
\[
\xymatrix{
\widehat{\mc A}\deco[\varpi]^\vee(\overline R_K) \ar[r] & \underline \omega_{\mc A/R} \otimes_R \overline R/\varpi \overline R \\
\mc D^{\cl}(\overline R_K) \ar@{^{(}->}[u] \ar[r] & \underline \omega_{(\mc D^{\cl})^\vee /R} \otimes_R \overline R / \varpi \overline R \ar@{^{(}->}[u]
}
\]
so we can study the map $\dlog_{\mc (D^{\cl})^\vee, \overline R}$. We now have that $\dlog_{\mc (D^{\cl})^\vee, \overline R}(y) \equiv \beta \bmod{\gamma}$, with $\val(\beta)=\frac{1-\val(\gamma)}{q-1}$, so $\dlog_{\mc (D^{\cl})^\vee, \overline R}(y)=0$ if and only if $\val(\gamma)\leq\frac{1-\val(\gamma)}{q-1}$, i.e.\ if and only if $\val(\gamma) \leq \frac{1}{q}$. If $y \in \mc C \deco(\overline R_K)^\perp \setminus \set{0}$, we have $\mc D^\perp = \mc C \deco (\overline R_K)$ and $\val(\gamma)=\val(E) \leq w \leq \frac{1}{q}$. It follows that $\mc C \deco(\overline R_K)^\perp$ is contained in the kernel of $\dlog_{1,\mc A}$. The first part of proposition follows since $\widehat{\mc A}\deco[\varpi]^\vee/({\mc C \deco(\overline R_K)}^\perp)^{\cl} \cong (\mc C \deco)^\vee$. If $y \not \in (\mc C \deco)(\overline R_K)^\perp$, the valuation of the points of $\mc D^\perp$ is $\frac{\val(E)}{q(q-1)}$. Looking at the map $\underline \omega_{(\mc D^{\cl})^\vee /R} \otimes_R \overline R / \varpi \overline R \hookrightarrow \underline \omega_{\mc A/R} \otimes_R \overline R/\varpi \overline R$, we see that, if $y \not \in (\mc C \deco)(\overline R_K)^\perp$, then $\dlog_{1,\mc A}(y) \equiv \beta \bmod{\varpi}$, with $\val(\beta)=\frac{q}{q-1}(1-\val(\gamma))=\frac{\val(E)}{q-1}\leq v$ as required. If $\widehat{\mc A}\deco$ has height $1$, a similar, but even simpler, argument, gives the result.
\end{proof}
\begin{rmk} \label{rmk: 1 over q}
In \cite{over}, the assumption $w < \frac{1}{p}$ is made from the very beginning to define the canonical subgroup. We need $w \leq \frac{1}{q}$ only to relate the map $\dlog$ with the canonical subgroup.
\end{rmk}
From now on, we will assume that $w \leq \frac{1}{q}$. We consider the morphism
\begin{gather*}
(\mc C \deco)_S = \Sp(S[x]/(x^q+\frac{\varpi}{E}x)) \to \mc{LT}_S= \Sp(S[[x]]) \\
E^{1/(q-1)} x \mapsfrom x
\end{gather*}
It gives a canonical non trivial point $\gamma \in (\mc C \deco)^\vee(S_K)$ that is a $\kappa$-generator of $(\mc C \deco)^\vee(S_K)$.
\begin{prop} \label{prop: gamma = E}
We have $\dlog_{1,S}(\gamma) \equiv {E_1}_{|\Spf(S)} \bmod{\varpi^{1-w}}$.
\end{prop}
\begin{proof}
Consider the following commutative diagram
\[
\xymatrix{
\mc (\mc A[\varpi]\deco)^\vee(S_K) \ar[r] \ar[d] & \underline \omega_{\mc A/R} \otimes_R S/\varpi^{1-w}S \ar[d] \\
(\mc C \deco)^\vee(S_K) \ar[r] \ar@{.>}[ur] & \underline \omega_{\mc C \deco/R} \otimes_R S/\varpi^{1-w}S
}
\]
Being $h \colon \mc C \deco \to \mc A[\varpi]\deco$ a closed immersion, the right vertical map is surjective. Both its domain and codomain are free $S/\varpi^{1-w} S$-module of rank $1$. It follows that the right vertical map is an isomorphism, so we can prove that $\dlog_{\mc C \deco,S}(\gamma) \equiv h^\ast({E_1}_{|\Spf(S)}) \bmod{\varpi^{1-w}}$. We have $\dlog_{\mc C \deco,S}(\gamma) =E^{1/(q-1)} \di(x)$, we conclude by Proposition~\ref{prop: descr diff}.
\end{proof}

\section{The Hodge-Tate sequence} \label{sec: HT}
We continue to work locally as in the previous section, using the same notations. We now need some results about $\widehat{\mc A}\deco[\varpi]^\vee$. Let $\underline \omega_{\mc A^\vee/R}$ be $\underline \omega_{\widehat{\mc A}\deco[\varpi]^\vee/R}$, and we choose $E' \in R$, a generator of this $R$-module.
\begin{prop} \label{prop: dual}
Let us suppose that $R$ is a discrete valuation ring, whose valuation extends the one of $\OP$. Then the valuation of $E'$ is the same as the valuation of $E$. Furthermore $(\mc C \deco(\overline R_K)^\perp)^{\cl}$ is the canonical subgroup of $\mc A\deco[\varpi]^\vee$.
\end{prop}
\begin{proof}
We can assume that $\widehat{\mc A}\deco$ has height $2$. We claim that the map $\dlog_{1,\mc A^\vee} \colon \mc A[\varpi]\deco(\overline R_K) \to \underline\omega_{\mc A^\vee/R} \otimes_R \overline R/\varpi \overline R$ has $\mc C \deco(\overline R_K)$ has kernel. Indeed, let $y \in \mc C \deco(\overline R_K)$, since the diagram
\[
\xymatrix{
\mc A\deco[\varpi](\overline R_K) \ar[r] & \underline \omega_{\mc A^\vee/R} \otimes_R \overline R/\varpi \overline R \\
\mc C \deco(\overline R_K) \ar@{^{(}->}[u] \ar[r] & \underline \omega_{(\mc C \deco)^\vee /R} \otimes_R \overline R / \varpi \overline R \ar@{^{(}->}
[u]
}
\]
is commutative, to prove that $\dlog_{1,\mc A^\vee}(y)=0$ it suffices to show that $\dlog_{(\mc C \deco)^\vee, \overline R}(y)=0$. But by \cite[Section~3]{faltings_strict}, we have $(\mc C \deco)^\vee \cong \Sp(R[x]/(x^q-Ex))$, so $\underline \omega_{(\mc C \deco)^\vee /R} \cong R/ER$. With this isomorphisms, we have $\dlog_{(\mc C \deco)^\vee, \overline R}(y) = \gamma$, with $\val(\gamma)=\frac{1-\val(E)}{q-1}\geq \val(E)$ since $\val(E) \leq \frac{1}{q}$. The claim follows since, by the analogue for $\varpi$-divisible groups of \cite[Lemma~2]{faltings_modular}, and Proposition~\ref{prop: ker dlog}, the kernel of $\dlog_{1,\mc A^\vee}$ is orthogonal to $\mc C \deco(\overline R_K)^\perp$ and hence has $\kappa$-dimension at most $1$. Using the analogue of the explicit calculations made in the proof of Proposition~\ref{prop: ker dlog}, we see that the fact that $\dlog_{1,\mc A^\vee}$ has a non trivial kernel implies that $\val(E') \leq \frac{1}{q}$. The statement about $(\mc C \deco(\overline R_K)^\perp)^{\cl}
$ follows. 
It remains to bound the valuation of $E'$, or equivalently, to bound the valuation of the points of $\mc C \deco(\overline R_K)^\perp$, that is $\frac{1-\val(E')}{q-1}$. Let us consider the isogeny $\widehat{\mc A}\deco[\varpi]^\vee \twoheadrightarrow \widehat{\mc A}\deco[\varpi]^\vee/(\mc C \deco(\overline R_K)^\perp)^{\cl}$. By \cite[Remark~2]{fargues_hodge}, it is given, after a suitable choice of coordinates, by the map
\begin{gather*}
R[[x]] \to R[[x]] \\
x \mapsto \prod_{\lambda \in \mc C \deco(\overline R_K)^\perp} (x - \lambda)
\end{gather*}
Since the valuation of the points of $\widehat{\mc A}\deco[\varpi]^\vee$ that are not in $\mc C \deco(\overline R_K)^\perp$ is $\frac{\val(E')}{q(q-1)}$, that is smaller than $\frac{1-\val(E')}{q-1}$, we have that the valuation of the image of these points under the isogeny is $\frac{\val(E')}{q-1}$. But $\widehat{\mc A}\deco[\varpi]^\vee/(\mc C \deco(\overline R_K)^\perp)^{\cl} \cong (\mc C\deco)^\vee$, whose points have valuation $\frac{\val(E)}{q-1}$, so $\val(E)=\val(E')$.
\end{proof}
\begin{rmk} \label{rmk: dual}
The above proposition implies that all our results about $\widehat{\mc A}\deco[\varpi]$ have an analogue for $\widehat{\mc A}\deco[\varpi]^\vee$, \emph{for the same constant $w$}.
\end{rmk}
We have the map
\[
\dlog_{\mc A} \colon \T_{\varpi}((\mc A[\varpi^\infty]\deco)^\vee)\otimes_{\OP} \widehat{\overline R} \to \underline \omega_{\mc A/R} \otimes_R \widehat{\overline R},
\]
and also its analogue for $(\mc A[\varpi^\infty]\deco)^\vee$,
\[
\dlog_{\mc A^\vee} \colon \T_{\varpi}(\mc A[\varpi^\infty]\deco)\otimes_{\OP} \widehat{\overline R} \to \underline \omega_{\mc A^\vee/R} \otimes_R \widehat{\overline R}.
\]
Let $\cdot^\ast$ mean `dual module', then we have an isomorphism of $\mc G$-modules
\[
\T_{\varpi}((\mc A[\varpi^\infty]\deco)^\vee)\cong \T_{\varpi}(\mc A[\varpi^\infty]\deco)^\ast(1),
\]
where $(\cdot)(1)$ means that the action of $\mc G$ is twisted by the Lubin-Tate character. We define $a_{\mc A}\colonequals  \dlog_{\mc A^\vee}^\ast(1)$.
\begin{defi} \label{defi: HT seq}
The \emph{Hodge-Tate sequence of $\mc A$} is the following sequence of $\widehat{\overline R}$-modules with semilinear action of $\mc G$:
\[
0 \to \underline \omega_{\mc A^\vee/R}^\ast \otimes_R \widehat{\overline R}(1) \stackrel{a_{\mc A}}{\longrightarrow} \T_{\varpi}((\mc A[\varpi^\infty]\deco)^\vee) \otimes_{\OP} \widehat{\overline R} \stackrel{\dlog_{\mc A}}{\longrightarrow} \underline \omega_{\mc A/R} \otimes_R \widehat{\overline R} \to 0.
\]
\end{defi}
\begin{prop} \label{prop HT compl}
The Hodge-Tate sequence of $\mc A$ is a complex.
\end{prop}
\begin{proof}
It is enough to show that $\Homol^0(\widehat{\overline R}(-1),\mc G)=0$. This follows by \cite[Proposition~3.1.8]{brinon} and \cite[Section~9]{faltings_strict}.
\end{proof}
\begin{lemma}[{(\cite[Proposition~2.0.3]{brinon})}] \label{lemma: 0 div}
We have that $\varpi$ is not a $0$-divisor in $\widehat{\overline R}$ and that the natural map $\overline R \to \widehat{\overline R}$ is injective.
\end{lemma}
\begin{prop} \label{prop: coker dlog}
The cokernel of the map $\dlog_{\mc A}$ is killed by $\varpi^v$ and $\im(\dlog_{\mc A})$ is a free $\widehat{\overline R}$-module of rank $1$. Furthermore, $\ker(\dlog_{\mc A})$ is a projective $\widehat{\overline R}$-module of rank $1$.
\end{prop}
\begin{proof}
Using Proposition~\ref{prop: ker dlog}, this is proved as the first step of the proof of \cite[Proposition~2.5]{over} (see also Lemma~2.5).
\end{proof}
\begin{lemma} \label{lemma: a inj}
The map $a_{\mc A}$ is injective.
\end{lemma}
\begin{proof}
By Remark~\ref{rmk: dual} and Proposition~\ref{prop: coker dlog}, we know that the cokernel of $\dlog_{\mc A^\vee}$ is killed by $\varpi^v$, so the same must be true for the kernel of $a_{\mc A}$, but by Lemma~\ref{lemma: 0 div} this implies that $\ker(a_{\mc A})=0$.
\end{proof}
Let $\mc D \deco$ be $\mc C \deco(\overline R_K)^\perp$. From now on we will omit $(\overline R_K)$ in the notation, it should be clear from the context whether we are talking about the group scheme or about the group of points. We also write $\overline R_z$ for $\overline R/\varpi^z \overline R$ (and similarly for other objects).
\begin{lemma} \label{lema: step 2}
We have a commutative diagram, with exact bottom row:
\[
\xymatrix@C=20pt{
& \underline \omega_{\mc A^\vee/R}^\ast \otimes_R \overline R_1(1) \ar[r] \ar[d] & \T_{\varpi}((\mc A[\varpi^\infty]\deco)^\vee) \otimes_{\kappa} \overline R_1 \ar[r] \ar@{=}[d] & \underline \omega_{\mc A/R} \otimes_R \overline R_1 & \\
0 \ar[r] & \mc D \deco \otimes_\kappa \overline R_1 \ar[r] & (\mc A[\varpi]\deco)^\vee \otimes_{\kappa} \overline R_1 \ar[r] & (\mc C\deco)^\vee \otimes_{\kappa} \overline R_1 \ar[u] \ar[r] & 0
}
\]
Furthermore we have an isomorphism $\ker(\dlog_{\mc A}) \cong \im(\dlog_{\mc A^\vee})^\ast(1)$.
\end{lemma}
\begin{proof}
Using Remark~\ref{rmk: dual}, this is proved as the second step of the proof of \cite[Proposition~2.4]{over}.
\end{proof}
\begin{teo} \label{teo: HT}
The homology of the Hodge-Tate sequence is killed by $\varpi^v$, and we have a commutative diagram of $\mc G$-modules, with exact rows and vertical isomorphisms,
\[
\xymatrix@C=15pt{
0 \ar[r] & \ker(\dlog_{\mc A})_{1-v} \ar[r] \ar[d]^\wr & \T_\varpi((\mc A[\varpi^\infty]\deco)^\vee) \otimes_{\OP} \overline R_{1-v} \ar[r] \ar@{=}[d] & \im(\dlog_{\mc A})_{1-v} \ar[r] \ar[d]^\wr & 0\\
0 \ar[r] & \mc D\deco \otimes_\kappa \overline R_{1-v} \ar[r] & (\mc A[\varpi]\deco)^\vee \otimes_\kappa \overline R_{1-v} \ar[r] & (\mc C \deco)^\vee \otimes_\kappa \overline R_{1-v} \ar[r] & 0
}
\]
Furthermore, $\ker(\dlog_{\mc A})$ is a free $\widehat{\overline R}$-module of rank $1$.
\end{teo}
\begin{proof}
Again this is proved using the same argument of the third step of the proof of \cite[Proposition~2.4]{over}, using Remark~\ref{rmk: dual}.
\end{proof}
\begin{prop} \label{prop: HT ord}
The Hodge-Tate sequence is exact if and only if $\mc A$ is ordinary.
\end{prop}
\begin{proof}
This follows from the calculations made in the proof of Proposition~\ref{prop: ker dlog}.
\end{proof}
Let $\mc H$ be $\Gal(\overline R_K/S_K)$. In the following, $\delta$ will be an element of $\overline R_1$ that satisfies $\dlog_{1,\mc A}(\gamma)= \delta \omega$ and $\tilde \delta \in \widehat{\overline R}$ will be a lifting of $\delta$. We can assume $\delta \in S/\varpi S$ and $\tilde \delta \in S$.
\begin{prop} \label{prop: HT S}
Let $\mc F(S) \subseteq \underline \omega_{\mc A/R}\otimes_R S$ be the submodule generated by $\tilde\delta\omega \otimes 1$.
\begin{enumerate}
 \item We have that $\mc F(S)$ is a free $S$-module of rank $1$, with basis $\tilde\delta\omega$ and $\mc F(S) \otimes_S \widehat{\overline R} \cong \im(\dlog_{A})$;
 \item the $S$-module $\im(\dlog_{\mc A})^{\mc H}$ is equal to $\mc F(S)$;
 \item there is an isomorphism $\mc F(S)_{1-v} \cong (\mc C \deco)^\vee \otimes_\kappa S_{1-v}$, its base change to $\widehat{\overline R}$ gives, via $\mc F(S) \otimes_S \widehat{\overline R} \cong \im(\dlog_{A})$, the isomorphism of Theorem~\ref{teo: HT};
 \item there is an isomorphism $\mc F(S)^\ast(1) \otimes_S \widehat{\overline R} \cong \ker(\dlog_{\mc A})$.
\end{enumerate}
Furthermore, all the above isomorphisms are $\mc H$-equivariant.
\end{prop}
\begin{proof}
This is proved as \cite[Proposition~2.6]{over}.
\end{proof}
\begin{lemma} \label{lemma: comp F}
Let $\Spf(R')$ be a small affine of $\formMoneH(w)$ and suppose that $R'$ is an $R$-algebra. We write $\mc A'$ for the base change of $\mc A$ to $R'$. Let $\Spf(S')$ be the inverse image of $\Spf(R')$ under the map $\formMoneHpi(w) \to \formMoneH(w)$, then we have a natural isomorphism $\mc F(S)\otimes_S S' \cong \mc F(S')$, compatible with $\underline \omega_{\mc A/R} \otimes_R R' \cong \underline \omega_{\mc A'/R'}$.
\end{lemma}
\begin{proof}
By functoriality of $\dlog$, we have a natural morphism $\im(\dlog_{\mc A})\otimes_{\widehat{\overline R}} \widehat{\overline{R'}} \to \im(\dlog_{\mc A'})$, that is compatible with the isomorphism $\underline \omega_{\mc A/R} \otimes_R R' \cong \underline \omega_{\mc A'/R'}$. Taking Galois invariants we obtain, by Proposition~\ref{prop: HT S}, a morphism $\mc F(S)\otimes_S S' \to \mc F(S')$, that is an isomorphism modulo $\varpi^{1-v}$ by Theorem~\ref{teo: HT}. The lemma follows.
\end{proof}
\section{\texorpdfstring{The sheaves $\Omega_w^\chi$ and $\Omega_{r,w}$}{The sheaves Omega w chi and Omega r,w}} \label{sec: sheaf}
We start assuming that $e \leq p-1$, we explain in Section~\ref{subsec: deep} how to remove this hypothesis.
\subsection{Generalities about locally analytic characters} \label{subsec: gen}
Let $A$ be a $K$-affinoid algebra. We will consider only $\FP$-locally analytic characters
\[
\chi \colon \OPast = \mu_{q-1} \times (1+\varpi\OP) \to A^\ast.
\]
\begin{defi} \label{defi: acc char}
Let $r \geq 1$ be an integer. We say that a character $\chi \colon \OP^\ast \to K^\ast$ is \emph{$r$-accessible} if it is of the form $t \mapsto [t]^i \langle t\rangle^s \colonequals [t]^i \exp(s\log(\langle t \rangle))$ for all $t$ with $\val(\langle t \rangle -1)\geq r$, where:
\begin{itemize}
 \item $i \in \Z/(q-1)\Z$;
 \item $[\cdot]$ is the Teichm\"uller character, and $[t]$ means $[\cdot]$ applied to the reduction of $t$ modulo $\varpi$;
 \item $\langle t\rangle\colonequals  t/[t]$ and $s \in K$ is such that $\val(s) > \frac{e}{p-1} - r$.
\end{itemize}
If $\chi$ is $1$-accessible, we will simply say that $\chi$ is accessible. In this case we write $\chi=(s,i)$. Given an integer $k$, we view it as the accessible character $t \mapsto t^k$. Note that any locally analytic character is $r$-accessible for some $r$.
\end{defi}
Let $\mc W$ be the \emph{weight space} for locally analytic characters: it is an $\FP$-rigid analytic space whose $A$-points, for any $\FP$-affinoid algebra $A$, are $\mc W(A) = \Hom_{\rm{loc-an}}(\OPast, A^\ast)$. There is a natural bijection between the set of connected components of $\mc W$ and $\Z/(q-1)\Z$. Let $\mc B$ be the component corresponding to the identity. We then have $\mc W = \coprod_{\Z/(q-1)\Z} \mc B$. By \cite[Theorem~3.6]{fourier}, we know that $\mc B$ is a twisted form, over $\C_p$, of $\mc B(1)$, the open disk of radius $1$. Note that $\mc B$ is isomorphic to $\mc B(1)$ if and only if $\FP=\Q_p$ (\cite[Lemma~3.9]{fourier}). In general $\mc B$ is a closed subvariety of $\mc B^d(1)$, the $d$-dimensional open polydisk of radius $1$, where $d =[\FP:\Q_p]$.

Let $t_1$ be $|\varpi|^{\frac{e}{p-1}}$, and, given an integer $r \geq 2$, we define $t_r < 1$ as the largest number such the following condition hold: if $x \in \C_p$ satisfies $|x-1| < t_r$, then $|\log(x)| < t_1 |\varpi|^{1-r}$. We have $t_r \to 1$ as $r \to \infty$. Let $\mc B^d(t_r)$ be the open $d$-dimensional polydisk of radius $t_r$. For $r \geq 1$, we fix $\mc D_r$, a closed ball such that $\mc B^d(t_{r-1}) \subset \mc D_r \subset \mc B^d(t_r)$. Viewing $\mc B$ as a subvariety of $\mc B^d(1)$, we define $\mc B_r$ via $\mc B_r \colonequals \mc B \cap \mc D_r$ (see \cite[Section~2]{fourier}). We write $\mc W_r$ for $\coprod_{\Z/(q-1)\Z} \mc B_r$.
\begin{lemma} \label{lemma: adm char} \label{lemma: v acc char}
Any $\chi \in \mc W_r(K)$ is an $r$-admissible character.
\end{lemma}
\begin{proof}
We may assume $\chi \in \mc B_r(K) \subseteq \mc B^d(t_r)(K)$. In this case we can take $s \colonequals \frac{\log(\chi(1+\varpi^r))}{\log(1+\varpi^r)}$.
\end{proof}
\begin{rmk}
We have that $\mc W_r$ is an affinoid subdomain of $\mc W$ and that $\set{\mc W_r}_{r \geq 0}$ is an admissible covering of $\mc W$. In particular, any character $\chi \in \mc W(K)$ lies in some $\mc W_r(K)$. Furthermore we know that any $\chi \in \mc W_r(K)$ is $r$-admissible.
\end{rmk}
\subsection{The case of accessible characters}
\begin{defi} \label{defi: sheaf F}
We write $\mc F$ for the unique locally free $\mc O_{\formMoneHpi(w)}$-module of rank $1$ that satisfies $\mc F(\Spf(S))=\mc F(S)$, for $\Spf(S)$ an open affine of $\formMoneHpi(w)$ as above. Here $\mc F(S)$ is the free $S$-module of rank $1$ defined in Section~\ref{sec: HT}.
\end{defi}
By Theorem~\ref{teo: HT}, we have an isomorphism of sheaves
\[
\mc F/\varpi^{1-v}\mc F \cong (\mathfrak C \deco)^\vee \otimes_\kappa \mc O_{\formMoneHpi(w)}/\varpi^{1-v}\mc O_{\formMoneHpi(w)}.
\]
\begin{defi} \label{defi: sheaf F'}
Using the above isomorphism, we define $\mc F'_v$ as the inverse image of the constant sheaf of sets $(\mathfrak C \deco)^\vee \setminus \set{0}$ under the natural map $\mc F \to \mc F/\varpi^{1-v}\mc F$.
\end{defi}
\begin{lemma} \label{lemma: F loc free}
Let $\Spf(S) \to \formMoneHpi(w)$ be an open affine, with associated abelian scheme $\mc A \to \Sp(S)$, and assume that $\underline \omega_{\mc A/S}$ is free, generated by $\omega$. Then we have that $\mc F(\Spf(S))$ is free, and $\omega^{\std}\colonequals {E_1}_{|\Spf(S)}$ gives a basis.
\end{lemma}
\begin{proof}
We can write ${E_1}_{|\Spf(S)}=E^{1/(q-1)}\omega$, for some $E\in S$. We can assume that $\Spf(S)$ is the inverse image of $\Spf(R) \to \formMoneH(w)$. Using the notations introduced before Proposition~\ref{prop: HT S}, we have that $\mc F(S)$ is generated by $\tilde \delta \omega$. By Proposition~\ref{prop: gamma = E}, we have $\tilde \delta \equiv E^{1/(q-1)} \bmod{\varpi^{1-w}}$. The lemma follows by completeness of $S$ and Nakayama's lemma.
\end{proof}
\begin{coro} \label{coro: F free}
We have that $\mc F$ is a \emph{free} $\mc O_{\formMoneHpi(w)}$-module of rank $1$, with $\omega^{\std}$ as basis.
\end{coro}
\begin{defi} \label{defi: sheaf S}
Let $\mc S_v$ be the sheaf of abelian groups, on $\formMoneHpi(w)$, defined by
\[
\mc S_v \colonequals  \OPast (1+\varpi^{1-v}\mc O_{\formMoneHpi(w)}).
\]
\end{defi}
\begin{prop} \label{prop: torsor}
We have that $\mc  F'_v$ is a Zariski $\mc S_v$-torsor and it is generated by $\omega^{\std}$.
\end{prop}
\begin{proof}
This follows from Proposition~\ref{prop: HT S} and Lemma~\ref{lemma: F loc free}.
\end{proof}
We write $\vartheta$ for the natural morphism $\vartheta \colon \formMoneHpi(w) \to \formMoneH(w)$. Its rigidification is Galois with $\kappa^\ast$ as Galois group. The action of $\kappa^\ast$ on $\vartheta^{\rig}$ extends to an action on $\vartheta$. Throughout this section, we fix an accessible character $\chi = (s,i)$. We will assume that
\[
w < (q-1)\left(\val(s) + 1 - \frac{e}{p-1} \right).
\]
Let $x$ be a local section of $\mc S_v$ over $\mathfrak V=\Spf(S)$. We can write $x=ub$, where $u$ is a section of $\OPast$ and $b$ is a section of $1+\varpi^{1-v}\mc O_{\mathfrak V}$. We set $x^\chi \colonequals \chi(u) b^s$, that is another section of $\mc S_v$. Note that, if $t \in 1+\varpi \OP$, we have $\chi(t) = t^s$, so $x^\chi$ is well defined. We will write $\mc O_{\formMoneHpi(w)}^{(\chi)}$ for $\mc O_{\formMoneHpi(w)}$ with the action of $\mc S_v$  by multiplication, twisted by $\chi$.

We have a natural action of $\mc S_v$ on $\mc F'_v$. In particular we can consider the sheaf
\[
\tilde \Omega_w^\chi \colonequals  \shom_{\mc S_v}(\mc F'_v, \mc O_{\formMoneHpi(w)}^{(\chi^{-1})}),
\]
where $\shom_{\mc S_v}(\cdot,\cdot)$ means homomorphisms of sheaves with an action of $\mc S_v$. By Proposition~\ref{prop: torsor}, we have that $\tilde \Omega_w^\chi$  is an invertible sheaf of $\mc O_{\formMoneHpi(w)}$-modules.
Note that, to specify $f$, a global section of $\tilde \Omega_w^\chi$, it is enough to give $f(\omega^{\std})$.
Since $\kappa^\ast$ acts on $(\mathfrak C \deco)^\vee \setminus \set{0}$, we have an action of $\kappa^\ast$ on $\mc F'_v$ and also an action of $\kappa^\ast$ on $\vartheta_\ast \mc O_{\formMoneHpi(w)}^{(\chi^{-1})}$.

Since $\vartheta$ is finite, we have that $\vartheta_\ast \tilde \Omega_w^\chi$ is a coherent sheaf of $\mc O_{\formMoneH(w)}$-modules. The action of $\kappa^\ast$ on $\mc F'_v$ and on $\vartheta_\ast \mc O_{\formMoneHpi(w)}^{(\chi^{-1})}$ gives an action of $\kappa^\ast$ on $\vartheta_\ast \tilde \Omega_w^\chi$. In particular, we have an action of $\kappa^\ast$ on the global section of $\tilde \Omega_w^\chi$. We will write this action by $f \mapsto f_{|\langle a \rangle}$, for $a \in \kappa^\ast$. These operators are the analogue of the usual diamond operators.
\begin{defi} \label{defi: Omega}
We define the sheaf $\Omega_w^\chi = \Omega_w^{(s,i)}$ on $\formMoneH(w)$ as $\Omega_w^\chi \colonequals  \left(\vartheta_\ast \tilde \Omega_w^\chi\right)^{\kappa^\ast}$.
\end{defi}
Let $\mathfrak V = \Spf(S) \to \formMoneHpi(w)$ be an open affine. We will write $X_{\chi,v}$ for the unique element of $\tilde \Omega_w^\chi(\mathfrak V)$ that satisfies $X_{\chi,v}(b \omega^{\std})=b^{-s}$, for all $b \in 1+\varpi^{1-v}S$. For various $\mathfrak V$'s, the $X_{\chi,v}$'s glue together, so we obtain a global section of $\tilde \Omega_w^\chi$, denoted again $X_{\chi,v}$.
\begin{lemma} \label{lemma: omega Hp loc free}
We have that $\tilde \Omega_w^\chi$ is a \emph{free} $\mc O_{\formMoneHpi(w)}$-module of rank $1$, with $X_{\chi,v}$ as basis.
\end{lemma}
\begin{proof}
This follows from Lemma~\ref{lemma: F loc free} and Proposition~\ref{prop: torsor}.
\end{proof}
\begin{rmk} \label{rmk: iso omega}
Let $\chi'=(s,j)$ be another accessible character (note that we have the same $s$ for $\chi$ and $\chi'$). We have a canonical isomorphism $\beta_{\chi,\chi'} \colon \tilde \Omega_w^\chi \stackrel{\sim}{\longrightarrow} \tilde \Omega_w^{\chi'}$, that sends $X_{\chi,v}$ to $X_{\chi',v}$. This isomorphism does not respect the action of $\kappa^\ast$, but we have that $\beta_{\chi,\chi'}$ induces an isomorphism $\tilde \Omega_w^\chi \cong \tilde \Omega_w^{\chi'}[j-i]$. Here by $\tilde \Omega_w^{\chi'}[j-i]$ we mean $\tilde \Omega_w^{\chi'}$ with the action of $\kappa^\ast$ twisted by $[\cdot]^{j-i}$.
\end{rmk}
\begin{defi} \label{defi: acc character Np}
We define the space of $\varpi$-adic modular forms with respect to $D$, level $\koneHpi$, weight $\chi$ and growth condition $w$, with coefficients in $K$, as
\[
S^D(K,w,\koneHpi,\chi)\colonequals \Homol^0(\formMoneHpi(w),\tilde \Omega_w^\chi)_K.
\]
If $\chi$ is an integer, by Lemma~\ref{lemma: int weight ok} below, we have $S^D(K,w,\koneHpi,\chi)=S^D(V,w,\koneHpi,\chi)_K$.
\end{defi}
\begin{prop} \label{prop: deco Omega}
There is a canonical $\kappa^\ast$-equivariant isomorphism of $\mc O_{\formMoneH(w)}$-modules
\[
\vartheta_\ast \tilde \Omega_w^\chi = \bigoplus_{j \in \Z/(q-1)\Z} \Omega_w^{(s,j)},
\]
such that $\Omega_w^{(s,j)}$ is the submodule of $\vartheta_\ast \tilde \Omega_w^\chi$ on which $\kappa^\ast$ acts via multiplication by $[\cdot]^{j-i}$.
\end{prop}
\begin{proof}
This is the analogue of \cite[Lemma~3.3]{over}. By Remark~\ref{rmk: iso omega}, $\Omega_w^{(s,j)}$ is equal to the set of invariant elements of $\vartheta_\ast \tilde\Omega_w^\chi[i-j]$, so it is the submodule of $\vartheta_\ast \tilde \Omega_w^\chi$ where $\kappa^\ast$ acts via $[\cdot]^{j-i}$. The order of $\kappa^\ast$ is $q-1$, that is invertible in all our rings, so $\vartheta_\ast \tilde \Omega_w^\chi$ can be decomposed, locally on $\formMoneH(w)$, as the direct sum of eigenspace of $\kappa^\ast$. The proposition follows.
\end{proof}
\begin{rmk} \label{rmk: indentification}
From now on we will use the above proposition to tacitly identify $\Omega_w^{(s,j)}$ with the submodule of $\vartheta_\ast \tilde \Omega_w^\chi$ on which $\kappa^\ast$ acts via $[\cdot]^{j-i}$.
\end{rmk}
\begin{coro} \label{coro: Omega invertible}
The rigidification of $\Omega_w^\chi$ is an invertible sheaf of $\mc O_{\formMoneH(w)^{\rig}}$-modules.
\end{coro}
\begin{defi} \label{defi: acc character}
We define the space of $\varpi$-adic modular forms with respect to $D$, level $\koneH$, weight $\chi$ and growth condition $w$, with coefficients in $K$, as
\[
S^D(K,w,\koneH,\chi)\colonequals \Homol^0(\formMoneH(w),\Omega_w^{(s,i)})_K.
\]
If $\chi$ is an integer, by Proposition~\ref{prop: int ok} below, we have $S^D(K,w,\koneH,\chi)=S^D(V,w,\koneH,\chi)_K$.
\end{defi}
Let $w' \geq w$ be a rational number that satisfies the same conditions of $w$. We set $v'\colonequals  \frac{w'}{q-1}$. We have natural morphisms $f_{w,w'} \colon \formMoneH(w) \to \formMoneH(w')$ and $g_{w,w'} \colon \formMoneHpi(w) \to \formMoneHpi(w')$.
\begin{lemma} \label{lemma: dir limi}
We have a natural isomorphism of $\mc O_{\formMoneHpi(w)}$-modules $\tilde \rho_{v,v'}\colon g^\ast_{w,w'}(\tilde \Omega_{w'}^\chi) \cong \tilde \Omega_w^\chi$. We have that $\tilde \rho_{v,v}= {\rm{id}}$ and, if $w'' \geq w'$ satisfies the same conditions of $w$, we have $\tilde \rho_{v,v''}=\tilde \rho_{v,v'}\circ g^\ast_{w,w'}(\tilde \rho_{v',v''})$, where $v''\colonequals \frac{w''}{q-1}$. Furthermore, we obtain a canonical morphism
\[
\rho_{v,v'}\colon f^\ast_{w,w'}(\Omega_{w'}^\chi) \to \Omega_w^\chi,
\]
that is an isomorphism after rigidification. The $\rho_{v,v'}$'s satisfy the same conditions as the $\tilde \rho_{v,v'}$'s do. Finally, we have $\rho_{v,v'}(X_{\chi, v'})=X_{\chi, v}$.
\end{lemma}
\begin{proof}
This is proved as \cite[Lemma~3.5]{over}.
\end{proof}
\begin{defi} \label{defi: overconvergent acc}
Using Lemma~\ref{lemma: dir limi}, we can define the space of \emph{overconvergent} modular forms with respect to $D$, level $\koneH$, weight $\chi$ and growth condition $w$, with coefficients in $K$, as
\[
S^D_\dagger(K,\koneH,\chi)\colonequals  \varinjlim_{w > 0} S^D(K,w,\koneH,\chi).
\]
\end{defi}
\begin{rmk} \label{rmk: funct Omega}
Let $i_{\mc A},i_{\mc B} \colon \Spf(S) \to \formMoneHpi(w)$ be two affine points of $\formMoneHpi(w)$. We write $\mc A$ and $\mc B$ for the abelian scheme corresponding to $i_{\mc A}$ and $i_{\mc B}$, respectively. Suppose we are given a morphism $f\colon\mc B \to \mc A$ over $S$. We obtain, by functoriality of $\dlog$, a morphism $\im(\dlog_{\mc A}) \to \im(\dlog_{\mc B})$ compatible with the natural pullback $\underline \omega_{\mc A/S} \to \underline \omega_{\mc B/S}$. Taking Galois invariants we obtain, by Proposition~\ref{prop: HT S}, a morphism $f^\ast \colon \mc F(i_{\mc A}(\Spf(S))) \to \mc F(i_{\mc B}(\Spf(S)))$. Let us now suppose that $f\colon \mc B \to \mc A$ is an isogeny, and that its kernel intersects trivially the canonical subgroup of $\mc B$. In this case we have a commutative diagram
\[
\xymatrix{
\mc F(i_{\mc A}(\Spf(S)))_{1-v} \ar[d]_\wr \ar[r]^{f^\ast_{1-v}} & \mc F(i_{\mc B}(\Spf(S)))_{1-v} \ar[d]^\wr \\
(\mc C (\mc A)\deco)^\vee \otimes_\kappa S_{1-v} \ar[r]^{(f\deco)^\vee} & \mc C ((\mc B)\deco)^\vee \otimes_\kappa S_{1-v}
}
\]
By assumption, $(f\deco)^\vee$ is an isomorphism, so $f^\ast$ is an isomorphism, modulo $\varpi^{1-v}$. This implies that $f^\ast$ is an isomorphism, so we have isomorphisms $\mc F'_v(i_{\mc A}(\Spf(S))) \cong \mc F'_v(i_{\mc B}(\Spf(S)))$ and
\[
\shom_{{\mc S_v}_{|i_{\mc B}(\Spf(S))}}({\mc F'_v}_{|i_{\mc B}(\Spf(S))},\mc O_{\Spf(S)}^{(\chi^{-1})}) \stackrel{\sim}{\to}
\shom_{{\mc S_v}_{|i_{\mc A}(\Spf(S))}}({\mc F'_v}_{|i_{\mc A}(\Spf(S))},\mc O_{\Spf(S)}^{(\chi^{-1})}).
\]
One can prove that $\shom_{{\mc S_v}_{|i_{\mc A}(\Spf(S))}}({\mc F'_v}_{|i_{\mc A}(\Spf(S))},\mc O_{\Spf(S)}^{(\chi^{-1})}) \cong i_{\mc A}^\ast \tilde \Omega_w^\chi$, and similarly for $\mc B$. In particular, we obtain an isomorphism $i_{\mc A}^\ast \tilde \Omega_w^\chi \to i_{\mc B}^\ast \tilde \Omega_w^\chi$. We will be more interested in its inverse
\[
\tilde f^\chi \colon i_{\mc B}^\ast \tilde \Omega_w^\chi \to i_{\mc A}^\ast \tilde \Omega_w^\chi.
\]
Let $\vartheta \colon \Spf(S) \to \Spf(R)$ be as above. We have a canonical isomorphism $\vartheta^{\rig,\ast} \Omega_w^\chi \cong \tilde \Omega_w^\chi$. In particular, we have the morphism
\[
f^\chi \colon  (\vartheta \circ i_{\mc B})^\ast \Omega_w^\chi \otimes_V K \to (\vartheta \circ i_{\mc A})^\ast \Omega_w^\chi \otimes_V K.
\]
In the case $\chi = (k,k)$ is an integer, via the isomorphism of Lemma~\ref{lemma: int weight ok}, $f^k $ is the pullback of the $k$-th power of the invariant differentials with respect to the isogeny.
\end{rmk}
\subsection{Modular forms of integral weight}
Let $k$ be an integer. If $\mathfrak V \subseteq \formMoneHpi(w)$, we write $\phi_{k,\mathfrak V} \colon (\tilde \Omega_w^k)(\mathfrak V) \to (\underline \omega_{\koneHpi}^{\otimes k})(\mathfrak V)$ be the map given by $\phi_{k,\mathfrak V}(f)=f(\omega^{\std})(\omega^{\std})^{\otimes k}$, for $f \in \tilde\Omega_w^k(\mathfrak V)$. We obtain a morphism $\phi_k \colon \tilde \Omega_w^k \to \underline \omega_{\koneHpi}^{\otimes k}$.
\begin{lemma} \label{lemma: int weight ok}
We have that $\phi_k \otimes_V K$ is an isomorphism.
\end{lemma}
\begin{proof}
The lemma follows since $\omega^{\std}\otimes 1$ is a generator of $\underline \omega_{\koneHpi} \otimes_V K$ by Theorem~\ref{teo: HT}.
\end{proof}
\begin{rmk} \label{rmk: coker phi}
By Proposition~\ref{prop: HT ord}, we see that $\phi_k$ is an isomorphism if and only if $w=0$. In general, by Theorem~\ref{teo: HT}, we have that $\coker(\phi_k)$ is killed by $\varpi^{kv}$.
\end{rmk}
\begin{lemma} \label{lemma: neben E}
We have that $\kappa^\ast$ acts on $E_1 \in S^D(K,w,\koneHpi,(1,1))$ via $[\cdot]^{-1}$.
\end{lemma}
\begin{proof}
Let $a \in \kappa^\ast$ and let $\Spf(S) \subseteq \formMoneHpi(w)$ be an open affine. We write $f$ for the element of $\tilde\Omega^{(1,1)}_w(\Spf(S))$ corresponding to ${E_1}_{|\Spf(S)}$. It is the morphism $f \colon (1+\varpi^{1-v}S){E_1}_{|\Spf(S)} \to \mc O_{\formMoneHpi(w)}^{(-1,-1)}(\Spf(S)) = S$ given by $f({E_1}_{|\Spf(S)})= 1$. By definition of $\omega^{\std}$, we see that $a^{-1}$ sends $\omega^{\std}$ to $[a]\omega^{\std}$. In particular, $f_{|\langle a \rangle}$ is the map that sends $\omega^{\std}$ to $a^\sharp(f([a]\omega^{\std}))$ (here $a=(a,a^\sharp)$ as morphism of ringed spaces). But $f \in \tilde\Omega^{(1,1)}_w(\Spf(S))$, so $f([a] \omega^{\std})= [a]^{-1}f(\omega^{\std})=[a]^{-1}$. So $f_{|\langle a \rangle}(\omega^{\std})= [a]^{-1}$ and $f_{|\langle a \rangle}= [a]^{-1}f$.
\end{proof}
We consider $\tilde \Omega_w^k$ as a subsheaf of $\underline \omega^{\otimes k}$, via $\phi_k$.
\begin{prop} \label{prop: int ok}
We have that $\underline \omega^{\otimes k, \rig}_{\koneH} = \Omega_w^{(k,k), \rig}$.
\end{prop}
\begin{proof}
We work locally, as in the proof of Lemma~\ref{lemma: neben E}. Let $f \otimes 1 \in \tilde \Omega_w^k(\Spf(S)) \otimes_V K$ be the element corresponding to $\omega^{\otimes k} \otimes 1 \in \underline \omega^{\otimes k}_R \otimes_V K$, where $\omega$ is a generator of $\underline \omega_R$, that we can assume to be free. Writing ${E_1}_{|\Spf(S)} = \frac{(-\varpi)^{1/(q-1)}}{\alpha}\omega$, we have that $f$ gives the map ${E_1}_{|\Spf(S)} \mapsto \left(\frac{\alpha}{(-\varpi)^{1/(q-1)}}\right )^k$. As in the proof of Lemma~\ref{lemma: neben E}, we have that
\[
f_{|\langle a \rangle}(\omega^{\std})= a^\sharp(f([a]\omega^{\std}))=[a]^{-k}a^\sharp\left(\frac{\alpha}{(-\varpi)^{1/(q-1)}}\right )^k.
\]
Furthermore, we have $a^\sharp(\alpha)= [a]\alpha$, so $f_{|\langle a \rangle}(\omega^{\std})=1$, hence $f_{|\langle a \rangle}=f$. This shows that $\underline \omega^{\otimes k}_{\koneHpi} \otimes_V K \subseteq \Omega_w^{(k,k)} \otimes_V K$. For the other inclusion, note that any element $f$ of $\tilde \Omega_w^k(\Spf(S)) \otimes_V K$ can be written as $f = s \omega$, for some $s \in S_K$. A calculation similar to the one above shows that, if $\kappa^\ast$ acts trivially on $f$, then $s \in R_K$ as required.
\end{proof}
\begin{rmk} \label{rmk: level N}
By Corollary~\ref{coro: Omega invertible}, we have a decomposition
\[
S^D(K,w,\koneHpi,\chi) = \bigoplus_{j \in \Z/(q-1)\Z} S^D(K,w,\koneH,(s,j)).
\]
Note that if $f$ is of level $\koneHpi$ and has integral weight, say $k$, we cannot identify it with a modular form of integral weight $k$ and level $\koneH$. Instead, $f$ will have components that are modular forms of level $\koneH$ and weight $x \mapsto \langle x \rangle ^k [x]^j$, for various $j \in \Z/(q-1)\Z$. This is very similar to the case of elliptic modular forms (see \cite[Sections I.3.4-7]{gouvea}). We have $X_{1,v}= E_1$, and we see that $X_{\chi,v}$ is a global section of $\Omega_w^{(s,0)}$. It follows that $X_{\chi,v}$, when considered as a modular form of level $\koneH$, has integral weight if and only if $s$ is an integer congruent to $0$ modulo $q-1$. For example we have that $X_{q-1,v}= E_{q-1}$ has weight $q-1$ as one expects.
\end{rmk}
\begin{rmk} \label{rmk: corr HT}
Fix an open affine $\mathfrak U = \Spf(R) \subseteq \formMoneH(w)$, and let $\mathfrak V = \Spf(S)$ be the inverse image of $\mathfrak U$ under $\vartheta$. We write $\mc A \to \Sp(S)$ for the corresponding abelian scheme. We have a $\kappa^\ast$-equivariant isomorphisms $\Omega_w^1 \cong \mc F$ and $\Omega_w^{-1} \cong \mc F^\ast$. In particular there is a `corrected' \emph{exact} Hodge-Tate sequence
\[
0 \to \Omega_w^{-1}(\mathfrak V) \otimes_R \widehat{\overline R}(1) \to \T_\varpi((\mc A[\varpi^\infty]\deco)^\vee) \otimes_\kappa \widehat{\overline R} \to \Omega_w^1(\mathfrak V) \otimes_R \widehat{\overline R} \to 0
\]
\end{rmk}
\subsection{Katz' modular forms}
We can describe our modular forms in a more familiar way, using `test objects'.
\begin{defi} \label{defi:test object}
A test object is a sextuple $(\mc A/S,i,\theta,\bar \alpha,Y,\eta)$, where:
\begin{itemize}
 \item $\Spf(S) \to \formMoneHpi(w)$ is an affine point, with $S$ a normal and $\varpi$-adically complete $V$-algebra;
 \item $(\mc A, i, \theta, \bar\alpha)$ is an object of the moduli problem of level $\koneHpi$, with $\mc A$ defined over $S$;
 \item $Y$ is a section of $\underline \omega_{\mc A/S}^{\otimes 1-q}$ that satisfies $YE_{q-1}=\varpi^w$;
 \item $\eta$ is a global section of the pullback of $\mc F'$ to $\Spf(S)$.
\end{itemize}
\end{defi}
\begin{prop} \label{prop: mod form rules}
To give an element $f$ of $S^D(K,w,\koneHpi,\chi)$ is equivalent to give a rule that assign to every test object $T=(\mc A/S,i,\theta,\bar \alpha,Y,\eta)$ an element $\tilde f(T) \in S_K$ such that:
\begin{itemize}
 \item $\tilde f(T)$ depends only on the isomorphism class of $T$;
 \item if $\varphi \colon S \to S'$ is a morphism of normal and $\varpi$-adically complete $V$-algebras, and we denote with $T'$ the base change of $T$ to $S'$, we have $\tilde f(T') = \varphi(\tilde f(T))$.
\end{itemize}
\end{prop}
\begin{proof}
This is proved as \cite[Lemma~3.10]{over}.
\end{proof}
\begin{coro} \label{coro: neben f tilde}
Let $f$ be in $S^D(K,w,\koneHpi,\chi)$. We have that $f \in S^D(K,w,\koneH,(s,j))$ if and only if, for any test object $T=(\mc A/S,i,\theta,\bar \alpha,Y,\eta)$, we have $\tilde f_{|\langle a \rangle}(T) = [a]^{j-i} \tilde f(T)$.
\end{coro}
\subsection{Canonical subgroups of higher rank and general characters}
In this section we fix an integer $r \geq 1$ and we suppose that $w < \frac{1}{q^{r-2}(q+1)}$.
\begin{prop} \label{prop: big can sub}
We have that $\formAoneH(w)[\varpi^r]$ has a canonical subgroup $\mathfrak C_r$ stable under the action of $D$. Furthermore $(\mathfrak C_r)\deco$ has order $q^r$ and $\mathfrak C_1 = \mathfrak C$.
\end{prop}
\begin{proof}
Let $\mc A \to \Sp(R)$ as above. We prove the proposition by induction on $r$, we already know the case $r=1$. By assumption, $\mc A[\varpi]$ admits a canonical subgroup $\mc C$. In \cite[Section~4.4 and Theorem~10.1]{pay_totally}, it is proved that $\mc A/\mc C$ is another object of the moduli problem, and that the $R$-point corresponding to it lies in $\formMoneH(qw)$. Since $qw \leq \frac{1}{q^{r-3}(q+1)}$, by induction hypothesis we have a canonical subgroup $\mc C'_{r-1} \subseteq \mc A/\mc C[\varpi^{r-1}]$. We define $\mc C_r$ to be the kernel of the composite map
\[
\mc A \to \mc A/\mc C \to (\mc A/\mc C)/\mc C'_{r-1}.
\]
\end{proof}
Using the canonical subgroups of higher rank, everything we have done for level $\koneHpi$ in Section~\ref{sec: level Hp} can be repeated for level $\koneHpin{r}$. In particular we have the rigid space $\formMoneHpin{r}(w)^{\rig}$, and the formal scheme $\formMoneHpin{r}(w)$.
\begin{prop} \label{prop: ker dlog_n}
Let $\Spf(R) \to \formMoneH(w)$ be as above. We have that the kernel of the map $\dlog_{r,\mc A} \colon (\mc A[\varpi^r]\deco)^\vee \to \underline \omega_{\mc A/R}\otimes_R \overline R_r$ is $\mc (\mc D_r)\deco \colonequals  ((\mc C_r)\deco)^\perp$.
\end{prop}
\begin{proof}
We prove the proposition by induction, the case $r=1$ follows by Proposition~\ref{prop: ker dlog}. We have that $\mc A[\varpi^r]\deco$ is an extension of $\mc A[\varpi]\deco$ and $\mc A[\varpi^{r-1}]\deco$, and the same is true for the canonical subgroups. The proposition follows from the functoriality of $\dlog$ and \cite[Corollary~1]{fargues_hodge}.
\end{proof}
\begin{prop} \label{prop: iso sheaf HT}
We have a natural $\mc G$-equivariant isomorphism
\[
\im(\dlog_{\mc A})_{r-v} \cong ((\mc C_r)\deco)^\vee \otimes_{\OP} \overline R_{r-v}.
\]
\end{prop}
\begin{proof}
Using Proposition~\ref{prop: ker dlog_n}, this is proved exactly in the same way as Theorem~\ref{teo: HT}.
\end{proof}
We have a natural morphism $\vartheta_r \colon \formMoneHpin{r}(w) \to \formMoneH(w)$. Its rigidification is Galois, with $G_r \colonequals (\OP/\varpi^r \OP)^\ast$ as Galois group. As above, we have that $G_r$ acts on $\vartheta_r$ too.

Let $\mathfrak U = \Spf(R) \subseteq \formMoneH(w)$ be an open affine. We will write $\mathfrak V_r = \Spf(S_r)$ for the inverse image of $\mathfrak U$ under $\vartheta_r$. We have that $(\mc C_r)\deco$ becomes constant over $S_{r,K}$. Furthermore there is a canonical point of $(\mc C_r)\deco$, defined over $S_r$. We can thus repeat what we have done for $\mc C \deco$, and we obtain an isomorphism of sheaves of $\mc O_{\formMoneHpin{r}(w)}$-modules
\[
\mc F /\varpi^{r-v} \mc F \cong ((\mathfrak C_r)\deco)^\vee \otimes_{\OP} \mc O_{\formMoneHpin{r}(w)}/\varpi^{r-v}\mc O_{\formMoneHpin{r}(w)}.
\]

We now fix $\set{\zeta_{n}}_{n \geq 1}$, a sequence of $\C_p$-points of $\mathcal{LT}$ such that the order of $\zeta_n$ is exactly $\varpi^n$. We assume that $\varpi\zeta_{n+1}=\zeta_n$ for each $r$, and that $\zeta_1$ is our fixed $(-\varpi)^{1/(q-1)}$. If $\zeta_{r} \in V$, we obtain $\gamma_r$, a canonical $S_r$-point of $((\mc C_r)\deco)^\vee$.

If $w$ is smaller than $1/(q^{r-2}(q+1))$, we define the sheaf $\mc F'_{r,v}$ on $\formMoneHpin{r}(w)$ as the inverse image of the constant sheaf of sets given by the subset of $((\mathfrak C_r)\deco)^\vee$ of points of order exactly $\varpi^r$. We have that $\mc F'_{r,v}$ is a Zariski $\mc S_{r,v}$-torsor, where $\mc S_{r,v}\colonequals \OP^\ast(1+\varpi^{r-v}\mc O_{\formMoneHpin{r}(w)})$.

We now fix $\chi$, an $r$-accessible character. We assume that $\zeta_r \in V$. Let $s$ be the element of $\C_p$ associated to $\chi$. We assume that $w$ is smaller than $1/(q^{r-2}(q+1))$, so we have the canonical subgroup of level $r$. Let $x=ub$ be a local section of $\mc S_{r,v}$. We have that $b^s\colonequals \exp(s\log(b))$ makes sense, so we can write $x^\chi \colonequals  \chi(u)b^s$. We write $\mc O_{\formMoneHpin{r}(w)}^{(\chi)}$ for the sheaf $\mc O_{\formMoneHpin{r}(w)}$ with the action of $\mc S_{r,v}$ twisted by $\chi$. We define the locally free sheaf of rank $1$
\[
\tilde \Omega_w^\chi \colonequals  \shom_{\mc S_{r,v}}(\mc F'_{r,v}, \mc O_{\formMoneHpin{r}(w)}^{(\chi^{-1})}).
\]
We have an action of $G_r$ on $\mc F'_{r,v}$ and on $\vartheta_{r,\ast} \mc O_{\formMoneHpin{r}(w)}^{\chi^{-1}}$. We obtain a coherent sheaves of $\mc O_{\formMoneH(w)}$-modules $\vartheta_{r,\ast} \tilde \Omega_w^\chi$. This sheaf is endowed with an action of $G_r$.
\begin{defi} \label{defi: Omega r}
We define the sheaf $\Omega_w^\chi$ on $\formMoneH(w)$ as $\Omega_w^\chi \colonequals  \left( \vartheta_{r,\ast} \tilde \Omega_w^\chi \right)^{G_r}$.
\end{defi}
Everything we have done in the case of an accessible character can be repeated for $\chi$. In particular we have modular forms, convergent and overconvergent, of weight $\chi$ and various levels.

Let $h$ be an integer with $r \geq h$. Suppose that $\chi$ is $h$-accessible. We can repeat the above construction starting with $\formMoneHpin{h}(w)$, obtaining another sheaf on $\formMoneHpi(w)$. For $r \geq h$, we consider the natural morphism $\vartheta_{r,h} \colon \formMoneHpin{r}(w) \to \formMoneHpin{h}(w)$. The rigidification of $\vartheta_{r,h}$ is Galois. Its Galois group is $G_{r,h} \subseteq G_r$, the image of $1+\varpi^{h}\OP$.
\begin{prop} \label{prop: ind r}
We have an isomorphism of $\mc O_{\formMoneH(w)} \otimes_V K$-modules
\begin{gather*}
\sigma_{r,h} \colon \left( \vartheta_{h,\ast} \shom_{\mc S_{h,v}}(\mc F'_{h,v}, \mc O_{\formMoneHpin{h}(w)}^{(\chi^{-1})})\otimes_V K \right)^{G_{h}} \cong \\
\cong \left( \vartheta_{r,\ast} \shom_{\mc S_{r,v}}(\mc F'_{r,v}, \mc O_{\formMoneHpin{r}(w)}^{(\chi^{-1})})\otimes_V K \right)^{G_{r}}. 
\end{gather*}
Furthermore $\sigma_{r,r}= {\rm{id}}$, and, if $t \leq h$ is an integer, we have $\sigma_{r,t}= \sigma_{h,t} \circ \sigma_{r,h}$.
\end{prop}
\begin{proof}
This is proved as \cite[Lemma~3.20]{over}.
\end{proof}
\subsection{\texorpdfstring{The sheaves $\Omega_{r,w}$}{The sheaves Omega r,w}}
We show that the sheaves $\Omega_w^\chi$ can be put in families. Let $\pi_i$, for $i=1,2$, be the natural projection from $\mc W_r \times \formMoneHpin{r}(w)^{\rig}$ to the $i$-th factor. We write $\mc S_{r,v}$ also for $\pi_2^{-1}(\mc S_{r,v})$ and $\mc F'_{r,w}$ also for $\pi_2^{-1}(\mc F'_{r,w})$. Let $x=ub$ be a section of $\mc S_{r,v}$. If $A\otimes B$ is a local section of $\mc O_{\mc W_r \times \formMoneHpin{r}(w)^{\rig}}$, we define $x(A\otimes B)$ to be the local section of $\mc O_{\mc W_r \times \formMoneHpin{r}(w)^{\rig}}$ that corresponds to the function $(\chi, z) \mapsto \chi(a)A(\chi)b^\chi B(z)$, for $\chi \in \mc W_r(T)$ and $z \in \formMoneHpin{r}(w)^{\rig}(T)$, where $T$ is any affinoid $K$-algebra. We define the sheaf
\[
\tilde \Omega_{r,w} \colonequals  \shom_{\mc S_{r,v}}(\mc F'_{r,v}, \mc O_{\mc W_r \times \formMoneHpin{r}(w)}).
\]
\begin{rmk} \label{rmk: desc omega_r}
It is possible to put also the $X_{\chi,v}$ in families. Let $\mathfrak V_{r} = \Spf(S_{r})$ be an open affine of $\formMoneHpin{r}(w)$ as always. We write $X_{r,v}$ for the element of $\tilde \Omega_{r,w}(\mc W_r \times \mathfrak V_{r}^{\rig})$ that satisfies $X_{r,v}(\omega^{\std})= 1$.
\end{rmk}
As in the case of a single character, we have that $G_{r}$ acts on $({\rm{id}}\times \vartheta_{r})_\ast \tilde \Omega_{r,w}$.
\begin{defi}
Let $r \geq 1$ be an integer, and let $w \leq 1/(q^{r-2}(q+1))$ be a rational number. On $\mc W_r \times \formMoneH(w)^{\rig}$, we define the sheaf $\Omega_{r,w}\colonequals  (({\rm{id}} \times \vartheta_{r})_\ast \tilde \Omega_{r,w})^{G_{r}}$.
\end{defi}
By construction we obtain the following
\begin{prop} \label{prop: prop fam}
The sheaves $\Omega_{r,w}$ are locally free sheaves of $\mc O_{\mc W_r \times \formMoneH(w)^{\rig}}$-modules of rank $1$. For any $\chi \in \mc W_r(K)$, we have a natural isomorphism
\[
(\chi,{\rm{id}})^\ast(\Omega_{r,w}) \cong \Omega_w^\chi.
\]
Furthermore, if $r_1$ and $r_2$ are integers greater than $0$ and $w_i \leq 1/(q^{r_i-2}(q+1))$, for $i=1,2$, are rational numbers, then the restrictions of $\Omega_{r_1,w_1}$ and $\Omega_{r_2,w_2}$ to $\mc W_{r_1} \cap \mc W_{r_2} \times \formMoneH(w_1)^{\rig} \cap \formMoneH(w_2)^{\rig}$ coincide.
\end{prop}
Any local section $f$ of $\Omega_{r,w}$ should be thought as a family of modular forms. Since $\formMoneH(w)^{\rig}$ is an affinoid, by Proposition~\ref{prop: prop fam} and Tate acyclicity Theorem, any modular form of weight $\chi$ lives in a $p$-adic family.
\begin{rmk} \label{rmk: funct fam}
Assume, as in Remark~\ref{rmk: funct Omega}, that we are given an isogeny $f \colon \mc B \to \mc A$, where $\mc A $ and $\mc B$ correspond to $i_{\mc A}, i_{\mc B} \colon \Spf(S_{r}) \to \formMoneHpin{r}(w)$. Suppose that the kernel of $f$ intersects trivially the canonical subgroup of $\mc B$. Writing $i_{\mc A}$ and $i_{\mc B}$ also for the maps $\mc W_r \times \Spf(S_{r})^{\rig} \to \mc W_r \times \formMoneHpin{r}(w)^{\rig}$, we obtain the families of morphisms $\tilde f_r \colon i_{\mc A}^\ast \tilde \Omega_{r,w} \to i_{\mc B}^\ast \tilde \Omega_{r,w}$ and $f_r \colon ((\operatorname{id} \times \vartheta_{r}^{\rig})\circ i_{\mc A})^\ast \Omega_{r,w} \to (((\operatorname{id} \times \vartheta_{r}^{\rig})\circ i_{\mc B})^\ast \Omega_{r,w}$.
\end{rmk}
\subsection{The deeply ramified case} \label{subsec: deep}
We now briefly explain what can be done without assuming that $e \leq p-1$. We have an isomorphism $\OP^\ast \cong \mu_{q-1} \times \mu_{p^n} \times \OP$ for some $n$. We can assume that $1+\varpi$ maps to $1$ under the maps $\OP^\ast \to \OP$ given by the above decomposition, so, with a little abuse of notation, we can write $\OP \cong (1+\varpi)^{\OP}$ (but note that the logarithm is not injective on $1+\varpi \OP$). In this way $\mc W$ becomes isomorphic to the disjoint union of $(q-1)p^n$ copies of $\mc B$ (see Section~\ref{subsec: gen}). We define the notion of $r$-accessible character as above, but only in the case $r \geq \frac{e}{p-1}$. In this way the definition of $\mc W_r$ can be adapted without problems. More importantly, if $\chi$ is $r$-accessible and $x$ is a local section of $\mc S_{r,v}$, we have that $x^s$ is a well defined section of $\mc S_{r,v}$. The rest of the theory goes smoothly. Thus, the real difference is 
that we do not have an integral structure for the space of modular forms of level $\koneHpin{r}$ and weight $\chi$ for any $r$, but only for $r$ big enough. However, if we invert $\varpi$ (i.e. if we take rigidification), the maps $\vartheta_r$ and $\vartheta_{r,h}$ are \'etale, furthermore we have a residual action of $G_r$ and $G_{r,h}$ on our sheaves, so there are no problems in this case.
\section{\texorpdfstring{The $\U$ operator}{The U operator}} \label{sec: U}
Let $\chi \colon \OP^\ast \to K^\ast$ be a character in $\mc W_r$ and let $0< w \leq 1/(q^{r-2}(q+1))$ be \emph{positive}.

Let $z$ be a point of $\formMoneHpin{r}(w)^{\rig}$, and let $L$ be its residue field, so $z$ comes from a morphism $\gamma_z \colon \Spm(L) \to \formMoneHpin{r}(w)^{\rig}$. We write $\tilde \gamma_z \colon \Spf(\mc O_L) \to \formMoneHpin{r}(w)$ for the rigid point associated to $z$. We have
\[
\Homol^0(\Spm(L),\gamma_z^\ast \tilde \Omega^\chi_w)=\Homol^0(\Spf(\mc O_L), \tilde \gamma_z^\ast \tilde \Omega^\chi_w) \otimes_{\mc O_L} L.
\]
We fix an identification $\Homol^0(\Spf(\mc O_L), \tilde \gamma_z^\ast \tilde \Omega^\chi_w) \cong \mc O_L$ and, if $f$ is an element of $\Homol^0(\Spm(L),\gamma_z^\ast \tilde \Omega^\chi_w)$, we define $|f|_z$ using the natural absolute value on $\mc O_L$. Let now $f$ be in $\Homol^0(\formMoneHpin{r}(w)^{\rig}, \tilde \Omega^\chi_w)$, we define $|f(z)|\colonequals |\gamma_z^\ast f|_z$, and we set
\[
|f|\colonequals \sup_{z \in \formMoneHpin{r}(w)^{\rig}}\set{|f(z)|}.
\]
\begin{prop} \label{prop: ON ban}
The $\sup$ defined above is always finite. In this way, $S^D(K,w,\koneHpin{r},\chi)$ becomes a potentially orthonormizable $K$-Banach module.
\end{prop}
\begin{proof}
Since $\formMoneHpin{r}(w)^{\rig}$ is an affinoid, the proposition follows by \cite[Lemma~2.14]{pay_curve}. 
\end{proof}
\begin{defi}
Let $M$ be a Banach $A$-module, where $A$ is an affinoid $K$-algebra. Following \cite[Part I, Section 2]{buzz_eigen}, we say that $M$ satisfies the \emph{property (Pr)}, if there is a Banach $A$-module $N$ such that $M \oplus N$ is potentially orthonormizable.
\end{defi}
\begin{coro} \label{coro: ortho}
The subspace $S^D(K,w,\koneH,\chi) \subseteq S^D(K,w,\koneHpin{r},\chi)$ satisfies property (Pr).
\end{coro}
To define the $\U$ operator we need to introduce another type of curves. We use the notations of Section~\ref{sec: shim curves}. We define
\[
\koneHpinq{r}\colonequals \set{\left(\begin{array}{cc} a & b \\ c & d\end{array}\right) \in \koneHpin{r} \mbox{ s.t. } b\equiv 0 \bmod{\varpi}}.
\]
In the case $K_{\mc P}=\koneHpinq{r}$, a choice of a level structure is equivalent to a choice of $(Q,D,\bar \alpha^{\mc P})$, where (here $(A,i,\theta,\alpha)$ is an object of the moduli problem for $\FP$-algebras):
\begin{enumerate}
 \item $Q$ is an $R$-point of exact $\OP$-order $\varpi^r$ in $A[\varpi^r]\deco$;
 \item $D$ is a finite and flat $\OP$-submodule of $A[\varpi^r]\deco$ of order $q$ which intersects the $\OP$-submodule scheme generated by $Q$ trivially.
\end{enumerate}
In this case, the curve $M_K$ will be denoted $\MoneHpinq{r}$, it is a proper and smooth scheme over $K$. There is a natural morphism $\pi_1 \colon \MoneHpinq{r} \to \MoneHpin{r}$, defined by the natural transformation of functors that forgets $D$. We have that $\pi_1$ is flat, and, since $\MoneHpin{r} \to \Sp(K)$ is proper, also $\pi_1$ must be proper. It follows that $\pi_1$ is finite.

Given $C$, a subgroup scheme of $A[q]$ of rank $q^{4N}$, stable under the action of $\OD$, we say, following \cite[Section~4.4]{pay_totally}, that it is of `type 2' if $C_2^2 \oplus \cdots \oplus C_m^2 = A[q]_2^2 \oplus \cdots \oplus A[q]_m^2$ and the isomorphism $\theta \colon A[q] \stackrel{\sim}{\longrightarrow} A[q]^{\Dual}$ sends $C \hookrightarrow A[q]$ to $(A[q]/C)^{\Dual}\hookrightarrow A[q]^{\Dual}$. Note that $C$, if it is of type 2, it is uniquely determined by $C\deco$. Given $D$, a finite and flat $\OP$-submodule of $A[\varpi]\deco$, we write $t_2(D)$ for the unique subgroup scheme of $A[q]$, of type 2, such that $t_2(D)\deco = D$. We can now define another morphism $\pi_2 \colon \MoneHpinq{r} \to \MoneHpin{r}$. At level of points, it is defined by taking the quotient over $t_2(D)$: in \cite[Section~4.4]{pay_totally}, it is shown how to put a level structure on $A/t_2(D)$, except for the point of exact $\OP$-order $\varpi^r$, but, since $D$ intersects trivially the $\OP$-submodule scheme 
generated by $Q$, we can take for it the image of $Q$ under the natural map $A \to A/t_2(D)$. We are interested in the analytification of $\pi_1$ and $\pi_2$, denoted respectively $\pi_1^{\rig}$ and $\pi_2^{\rig}$.

The rigid space associated to $\MoneHpinq{r}$ will be denoted $\formMoneHpinq{r}^{\rig}$. Furthermore, we write $\formMoneHpinq{r}(w)^{\rig}$ for $(\pi_1^{\rig})^{-1}(\formMoneHpin{r}(w)^{\rig})$. We define the formal model $\formMoneHpinq{r}(w)$ as the normalization, via $\pi_1^{\rig}$, of $\formMoneHpin{r}(w)$ in $\formMoneHpinq{r}(w)^{\rig}$. In this way we obtain a formal model of $\pi_1^{\rig}$ , denoted $\mathfrak p_1 \colon \formMoneHpinq{r}(w) \to \formMoneHpin{r}(w)$.
\begin{prop} \label{prop: mod MoneHpiq}
Let $R$ be a normal and $\varpi$-adically complete $V$-algebra. Then there is a natural bijection between $\formMoneHpinq{r}(w)(R)$ and the set of isomorphism classes of sextuples $(\mc A, i, \theta, \bar\alpha, Y, \mc D)$, where:
\begin{itemize}
 \item $(\mc A, i, \theta, \bar\alpha, Y)$ is an object of the moduli problem, with $\mc A$ defined over $R$, of $\intMoneHpin{r}(w)$;
 \item $\mc D$ is a finite and flat $\OP$-submodule of $\mc A[\varpi^r]\deco$ of rank $r$ that intersects trivially the canonical subgroup of $\mc A[\varpi^r]\deco$.
\end{itemize}
\end{prop}
\begin{proof}
This can be proved as \cite[Lemma~3.11]{over}.
\end{proof}
\begin{lemma} \label{lemma: comp p1 p2}
Let $R$ be a normal and $\varpi$-adically complete $V$-algebra and let $(\mc A, i, \theta, \bar\alpha, Y, \mc D)$ be in $\formMoneHpinq{r}(qw)(R)$. Taking the quotient over $t_2(\mc D)$, we obtain an object of $\formMoneHpin{r}(w)(R)$.
\end{lemma}
\begin{proof}
It is enough to consider the case $r=1$. We can assume that $R$ is a a discrete valuation ring, whose valuation extends the one of $V$ and that $\mc A$ is supersingular. Let $\mc B$ be $\mc A/t_2(\mc D)$. Forgetting the extra structure, we need to prove that the $R$-point corresponding to $\mc B$ lies in $\formMoneH(qw)$. To prove this, let us consider the commutative diagram
\[
\xymatrix{
(\mc A[\varpi]\deco)^\vee (\overline R_K) \ar[r] & \underline \omega_{\mc A/R} \otimes_R \overline R_1 \ar[r]^(.65)\sim & \overline R_1 \\
(\mc B[\varpi]\deco)^\vee (\overline R_K) \ar[r] \ar@{->>}[u] & \underline \omega_{\mc B/R} \otimes_R \overline R_1 \ar[r]^(.65)\sim \ar[u] & \overline R_1 \ar[u]
}
\]
We use the notation of the proof of Proposition~\ref{prop: ker dlog}. We need to prove that $\val(E_{\mc B})\leq \frac{\val(E_{\mc A})}{q}$. The right vertical map is the reduction of the multiplication by an element of valuation $\frac{\val(E_{\mc A})}{q}$ by \cite[Remark~2]{fargues_hodge}. Looking at the proof of Proposition~\ref{prop: ker dlog}, we see that the image of the compositions of the horizontal maps are generated by elements of valuation, respectively, $\frac{\val(E_{\mc A})}{q-1}$ and $\frac{\val(E_{\mc B})}{q-1}$, so $\frac{\val(E_{\mc B})}{q-1} + \frac{\val(E_{\mc A})}{q} =\frac{\val(E_{\mc A})}{q-1}$ as required.
\end{proof}
Taking the quotient over $\mc D$, we define, on points, the morphism
\[
\mathfrak p_2 \colon \formMoneHpinq{r}(qw) \to \formMoneHpin{r}(w),
\]
Let $\formAoneHpiqn{r}(w)$ be the base change, via $\mathfrak p_1$, to $\formMoneHpinq{r}(w)$, of $\formAoneHpin{r}(w)$. We have that $\formAoneHpiqn{r}(w)$ is equipped with $\mathfrak D$, a subgroup of order $q$ of its $\varpi^r$-torsion, that intersects trivially its canonical subgroup. The isogeny
\[
\pi_{\mathfrak D} \colon \formAoneHpiqn{r}(qw) \to \formAoneHpiqn{r}(qw)/\mathfrak D
\]
is defined over $\formMoneHpinq{r}(qw)$. Since $\formAoneHpiqn{r}(qw)/\mathfrak D$ is the base change, via $f_{w,qw} \circ \mathfrak p_2$, to $\formMoneHpinq{r}(qw)$, of $\formAoneHpin{r}(qw)$, we obtain, by Remark~\ref{rmk: funct Omega} and Lemma~\ref{lemma: dir limi}, a morphism
\[
\tilde \pi_{\mathfrak D}^\chi \colon \mathfrak p_2^\ast \tilde \Omega^\chi_w \to \mathfrak p_1^\ast \tilde \Omega^\chi_{qw}.
\]
Let us consider $\tilde \U$, defined as the composition
\begin{gather*}
\Homol^0(\formMoneHpin{r}(qw), \tilde \Omega_{qw}^\chi \otimes_V K) \stackrel{\tilde \rho_{qw,w}^{\rig}}{\longrightarrow} \Homol^0(\formMoneHpin{r}(w), \tilde \Omega_w^\chi \otimes_V K) \stackrel{\mathfrak p_2^\ast}{\longrightarrow} \\
\to \Homol^0(\formMoneHpinq{r}(qw),\mathfrak p_2^\ast \tilde \Omega^\chi_w \otimes_V K) \stackrel{\tilde \pi_{\mathfrak D}^\chi}{\longrightarrow} \\
\to \Homol^0(\formMoneHpinq{r}(w),\mathfrak p_1^\ast \tilde \Omega^\chi_{qw}\otimes_V K) \stackrel{\pi_{1,\ast}^{\rig}}{\longrightarrow} \Homol^0(\formMoneHpin{r}(qw),\tilde \Omega_{qw}^\chi \otimes_V K),
\end{gather*}
where $\pi_{1,\ast}^{\rig}$ is the map induced by the trace, that is well defined since $\pi_1^{\rig}$ is finite and flat. All the maps used to define $\tilde \U$ are $G_{r}$-equivariant, so the same holds for $\tilde \U$. Taking $G_{r}$-invariants we obtain, from $\tilde \U$, a map $S^D(K,w,\koneH,\chi) \to S^D(K,w,\koneH,\chi)$.
\begin{defi} \label{defi: U}
Let $\chi$ be an $r$-accessible character. The map
\[
\U \colon S^D(K,qw,\koneH,\chi) \to S^D(K,qw,\koneH,\chi)
\]
is defined as $1/q$ times the map induced by $\tilde \U$.
\end{defi}
\begin{prop} \label{prop: U compl cont}
The operator $\U$ is completely continuous.
\end{prop}
\begin{proof}
We claim that $\tilde \U$ is completely continuous. Since $\tilde \U$ factors through $\tilde \rho_{qw,w}^{\rig}$, it is enough to prove that $\tilde \rho_{qw,w}^{\rig}$ is completely continuous, and this can be done in exactly the same way as \cite[Proposition~2.20]{pay_curve}. The proposition follows.
\end{proof}
\begin{rmk} \label{rmk: heck test}
Let us suppose that $r=1$. Let $f$ be an element of $S^D(K,w,\koneH,\chi)$. Take any test object $T=(\mc A/S, i, \theta, \bar\alpha,Y,\eta)$ as in Proposition~\ref{prop: mod form rules}. Let $S'$ be a normal and $\varpi$-adically complete $S$-algebra such that $S_K \to S'_K$ is finite and \'etale and all finite and flat subgroup schemes of $\mc A_{\overline S,K}[\varpi]\deco$ are defined over $S'_K$. Repeating what we have done in the proof of Proposition~\ref{prop: mod MoneHpiq}, we see that any finite and flat subgroup scheme of $\mc A_{S',K}[\varpi]\deco$ extends to a subgroup scheme of $\mc A_{S'}[\varpi]\deco$. Let $\mc D$ be any such subgroup, and suppose that $\mc D$ intersects trivially the canonical subgroup of $\mc A_{S'}[\varpi]\deco$. We have that $T$ gives a test object $((\mc A_{S'} / t_2(\mc D))/S', i', \theta', \bar\alpha',Y',\eta')$. Indeed the only non trivial thing to define is $\eta'$. Let $i_1,i_2 \colon \Spf(S) \to \formMoneHpi(w)$ be the morphisms corresponding to $\mc A$ and $\mc A/t_
2(\mc D)$, respectively. In Remark~\ref{rmk: funct Omega} we showed that there is an isomorphism between the global sections of $i_1^\ast \mc F'$ and $i_2^\ast \mc F'$. We define $\eta'$ as the image of $\eta$ under this isomorphism. We have
\[
\widetilde {f_{|\U}}(T) = \frac{1}{q} \sum_{\mc D} \tilde f(((\mc A_{S'} / t_2(\mc D))/S', i', \theta', \bar\alpha',Y',\eta')).
\]
\end{rmk}
For various $w$'s, the norms defined on $S^D(K,w,\koneH,\chi)$ are compatible, so $S^D_\dagger(K,\koneH,\chi)$ is naturally a Fr\'echet space, and we obtain a continuous operator $\U$ on it.

Using the maps $\tilde \pi_{\mathfrak D,r}$ defined in Remark~\ref{rmk: funct fam}, we can work with families: for any integer $r \geq 1$, we obtain an operator
\[
\tilde \U_r \colon \tilde \Omega_{r,w}(\mc W_r \times \formMoneHpin{r}(w)^{\rig}) \to \tilde \Omega_{r,w}(\mc W_r \times \formMoneHpin{r}(w)^{\rig}),
\]
such that the pullback $(\chi, {\rm{id}})^\ast (\tilde \U_r)$, for $\chi \in \mc W_r(K)$, is the $\tilde \U$ operator defined above. Everything we did above can be repeated for families, in particular we have the $\U_r$ operator and the following
\begin{prop} \label{prop: fam op compl cont}
For any integer $r \geq 1$ and any rational $w \leq 1/(q^{r-2}(q+1))$, we have that $\Homol^0(\Omega_{r,w}, \mc W_r \times \formMoneH(w)^{\rig})$ is a Banach $\mc O_{\mc W_r}(\mc W_r)$-module that satisfies the property (Pr). Furthermore the $\U_r$ operator is completely continuous.
\end{prop}
Kassaei has proved a result of classicality for modular forms of level $\koneH$ and integral weight $k$. Let $f$ be in $S^D(K,w,\koneH,k)$ and suppose that $f_{|\U}= a f$, for some $a \in K$. If $a$ satisfies $\val(a) < k - ef$, then $f$ is classical, i.e.\ it can be extended to a global section of $\underline \omega^{\otimes k}_{\formMHpi^{\rig}}$. See \cite[Theorem~5.1]{pay_curve}.

Let $\chi \colon \OP^\ast \to K^\ast$ be a locally analytic character and let $\nu \in \R$. Let $\mc V = \Spm(R) \subseteq \mc W$ be an affinoid that contains $\chi$. Using the notations of \cite[page~31]{bellaiche}, we have that $\Homol^0(\mc V \times \formMoneH(w)^{\rig}, \Omega_{r,w})^{\leq \nu}$ makes sense if $\mc V$ is sufficiently small. We have an isomorphism
\[
\Homol^0(\mc V \times \formMoneH(w)^{\rig}, \Omega_{r,w})^{\leq \nu} \otimes_{R,\chi} K \cong S^D(K,w,\koneH, \chi)^{\leq \nu}.
\]
In particular, we have the following
\begin{prop} \label{prop: defo}
Let $\nu$ be in $\R$ and let $f$ be in $S^D(K,w,\koneH, \chi)^{\leq v}$. Then there is an affinoid $\mc V \subseteq \mc W$ such that $f$ can be deformed to a family of modular forms over $\mc V$. Furthermore, the $\U$-operator acts with slope $\leq \nu$ on this family.
\end{prop}
\subsection{Other Hecke operators}
We now sketch the definition of other Hecke operators. Let $l \neq p$ be a rational prime. We write $\mc L_1, \ldots, \mc L_k$ for the primes of $F$ lying over $l$. Let $\mc L$ be $\mc L_1$. We assume that $l$ splits in $\Q(\sqrt{\lambda})$, and that $B$ is split at $\mc L$. We denote the completion of $F$ at $\mc L_i$ with $F_{\mc L_i}$. We have
\[
G(\Q_l) \cong \Q_l^\ast \times \GL_2(F_{\mc L}) \times \GL_2(F_{\mc L_2}) \times \cdots \times \GL_2(F_{\mc L_k}).
\]
In this section we make the assumption that the compact open subgroup $H$ is of the form
\[
H=\Z_l^\ast \times \GL_2(\mc O_{F_{\mc L}}) \times H'.
\]
Let $\varpi_l$ be a uniformizer of $\mc O_{F_{\mc L}}$. If $A$ is an abelian scheme as above, we have a decomposition of $A[\varpi_l]$ similar to that of $A[\varpi]$, so $A[\varpi_l]\deco$ is defined and it has an action of $\kappa_l \colonequals \mc O_{F_{\mc L}}/\varpi_l$.

Let $\chi \colon \OP^\ast \to K^\ast$ be an $r$-accessible character. Let $H_{\mc L}$ be the set of invertible $2\times2$ matrices with left lower corner congruent to $0$ modulo $\varpi_l$. The Shimura curve corresponding to the case $K_{\mc P}=\koneHpin{r}$ and $H = \Z_l^\ast \times H_{\mc L} \times H'$ will be denoted with $X$. We have that $X$ parametrizes objects of the moduli problem of $\MoneHpin{r}$ plus a finite and flat subgroup of $A[\varpi_l]\deco$ of order $|\kappa_l|$, stable under the action of $\mc O_{F_{\mc L}}$. If $D$ is such a subgroup, we can define $t_2(D)$ as in the case of subgroups of $A[\varpi]\deco$, and also the quotient of $A$ by $t_2(D)$ can be defined as in \cite[Section~4.4]{pay_totally}. We can repeat everything we have done for the $U$ operator. In particular we obtain $\mathfrak p_1, \mathfrak p_2 \colon \mathfrak X(w) \to \formMoneHpin{r}(w)$. Furthermore, we have a morphism $\tilde \pi_{\mathfrak D} \colon \mathfrak p_2^\ast \tilde \Omega_w^\chi \to \mathfrak p_1^\ast \tilde \Omega_w^\chi$.
\begin{defi} \label{defi: T_l}
We define the operator
\[
{\T}_{\mc L} \colon S^D(K,w,\koneHpin{r},\chi) \to S^D(K,w,\koneHpin{r},\chi)
\]
exactly as in the case of $\U$ (using $|\kappa_l| + 1$ as normalization factor).
\end{defi}
\begin{rmk} \label{rmk: T_l}
Note that $\tilde {\T}_{\mc L}$ is a continuous operator that it is not completely continuous.

Also the operators $\tilde {\T}_{\mc L}$ can be put in families. Furthermore, if $\chi$ is accessible, we have a description of $\tilde {\T}_{\mc L}$ in terms of testing objects similar to that of Remark~\ref{rmk: heck test}.
\end{rmk}

Let $r \geq 1$ be an integer, and assume that $0 < w$ is a rational number sufficiently small. Let $\mc Z_r$ be the spectral variety associated to the $\U$-operator acting on $\Homol^0(\mc W_r \times \formMoneH(w)^{\rig}, \Omega_{r,w})$. We have proved that all assumptions needed to use the machine developed by Buzzard in \cite{buzz_eigen} are satisfied, so we have the following
\begin{teo} \label{teo: eigen}
We have a rigid space $\mc C_r \subseteq \mc W_r \times \m A^{1,\rig}_K$ equipped with a finite morphism $\mc C_r \to \mc Z_r$. If $L$ is a finite extension of $K$, then the points of $\mc C_r(L)$ correspond to systems of eigenvalues of modular forms with growth condition $w$ and coefficients in $L$. If $x \in \mc C_r(L)$, let $\mc M(w)_x$ be the set of modular forms corresponding to $x$. Then all the elements of $\mc M(w)_x$ have weight $\pi_1(x) \in \mc W(L)$ and the $\U$-operator acts on $\mc M(w)_x$ with eigenvalue $\pi_2(x)^{-1}$. For various $r$ and $w$, these construction are compatible. Letting $r \to \infty$ we have $w \to 0$ and we obtain the global eigencurve $\mc C \subseteq \mc W \times \m A^{1,\rig}_K$.
\end{teo}

\bibliographystyle{amsalpha}
\bibliography{biblio}

\end{document}